\documentclass[12pt,sort&compress]{elsarticle}
\usepackage[utf8]{inputenc}

\usepackage{bm}

\usepackage{mathrsfs}
\usepackage{amsmath}
\usepackage{amsfonts}
\usepackage{amsthm}
\usepackage{color}
\usepackage{sansmath}

\usepackage{array}

\usepackage{caption}
\usepackage{subcaption}

\usepackage{float}
\usepackage{url}
\usepackage{multicol}
\usepackage[top=1in,bottom=1in,left=1in,right=1in]{geometry}
\usepackage[pdfborder={0 0 0},colorlinks,allcolors=blue]{hyperref}
\usepackage{txfonts}
\usepackage{setspace}
\usepackage{arydshln}
\usepackage{enumitem}
\usepackage{lipsum}
\usepackage{siunitx,booktabs}
\usepackage{nth}
\usepackage{accents} 

\theoremstyle{definition}
\newtheorem{remark}{Remark}

\allowdisplaybreaks

\providecommand{\e}[1]{\ensuremath{\times 10^{#1}}}

\newcommand{\etal}{\textit{et al}. }

\usepackage{prettyref}
\newrefformat{fig}{Fig.~\ref{#1}}
\newrefformat{Fig}{Fig.~\ref{#1}}
\newrefformat{Eq}{Eq.(\ref{#1})}
\newrefformat{eq}{Eq.(\ref{#1})}
\newrefformat{sec}{Section~\ref{#1}}
\newrefformat{tab}{Table~\ref{#1}}
\newrefformat{Tab}{Table~\ref{#1}}
\newrefformat{alg}{Algorithm~\ref{#1}}
\newrefformat{Alg}{Algorithm~\ref{#1}}

\usepackage[dvipsnames]{xcolor}
\newcommand\hl[1]{%
	\bgroup
	\hskip0pt\color{red!80!black}%
	#1%
	\egroup
}

\newtheorem{definition}{Definition}

\usepackage{algorithm}
\usepackage{algpseudocode}
\usepackage{algorithmicx}
\definecolor{shadecolor}{cmyk}{0,0,0,0.03}

\newcommand*\pnl{\Statex$\vartriangleright$ }
\newcommand*\qpnl{\Statex\quad$\vartriangleright$ }
\newcommand*\qqpnl{\Statex\qquad$\vartriangleright$ }

\newcommand*\subr[1]{\Statex\textbf{#1}}
\newcommand*\qcomnt[1]{\Statex \textit{#1}}
\newcommand*\qqcomnt[1]{\Statex\quad~~~\textit{#1}}
\renewcommand{\For}[3]{\Statex \textbf{for} \textit{#1} $\leftarrow$ \textit{#2} \textbf{to} \textit{#3} \textbf{do}}
\renewcommand*\EndFor{\Statex \textbf{end for}}

\newcommand*\qFor[3]{\Statex\quad\textbf{for} \textit{#1} $\leftarrow$ \textit{#2} \textbf{to} \textit{#3} \textbf{do}}
\newcommand*\EndqFor{\Statex\quad \textbf{end for}}

\newcommand*\qIf[3]{\Statex\quad\textbf{if} $#1 #2 #3$ \textbf{then}}

\newcommand*\EndqIf{\Statex\quad\textbf{end if}}

\journal{Computer Methods in Applied Mechanics and Engineering}

\begin{document}

\begin{frontmatter}

\title{A smoothed particle hydrodynamics approach for phase field modeling of brittle fracture}

\author[sbu]{Mohammad~Naqib~Rahimi}

\author[sbu,iacs]{Georgios~Moutsanidis\corref{cor1}}
\ead{georgios.moutsanidis@stonybrook.edu}
\cortext[cor1]{Corresponding author}

\address[sbu]{Department of Civil Engineering, Stony Brook University, Stony Brook, NY 11794, USA}

\address[iacs]{Institute for Advanced Computational Science, Stony Brook, NY 11794, USA}

\begin{abstract}
Fracture is a very challenging and complicated problem with various applications in engineering and physics. Although it has been extensively studied within the context of mesh-based numerical techniques, such as the finite element method (FEM), the research activity within the Smoothed Particle Hydrodynamics (SPH) community remains scarce. SPH is a particle-based numerical method used to discretize equations of continuum media. Its meshfree nature makes it ideal to simulate fracture scenarios that involve extreme deformations. However, to model fracture, SPH researchers have mostly relied  on ad-hoc empirical local damage models, cohesive zone approaches, or pseudo-spring models, which come with a set of drawbacks and limitations. On the other hand, phase field models of brittle fracture have recently gained popularity in academic circles and provide significant improvements compared to previous approaches. These improvements include the derivation from fundamental fracture theories, the introduction of non-locality, and the ability to model multiple crack initiation, propagation, branching, and coalescence, in situations where no prior knowledge of the crack paths is available. Nevertheless, phase field for fracture has not been studied within SPH. In this proof-of-concept paper we develop and implement a phase field model of brittle fracture within the context of SPH. Comprehensive mathematical and implementation details are provided, and several challenging numerical examples are computed and illustrate the proposed method's ability to accurately and efficiently simulate complex fracture scenarios.
\end{abstract}

\begin{keyword}
SPH; Particle methods; Phase field; Fracture mechanics; Non-local methods; Contact mechanics
\end{keyword}

\begin{graphicalabstract}
\includegraphics[width=\linewidth]{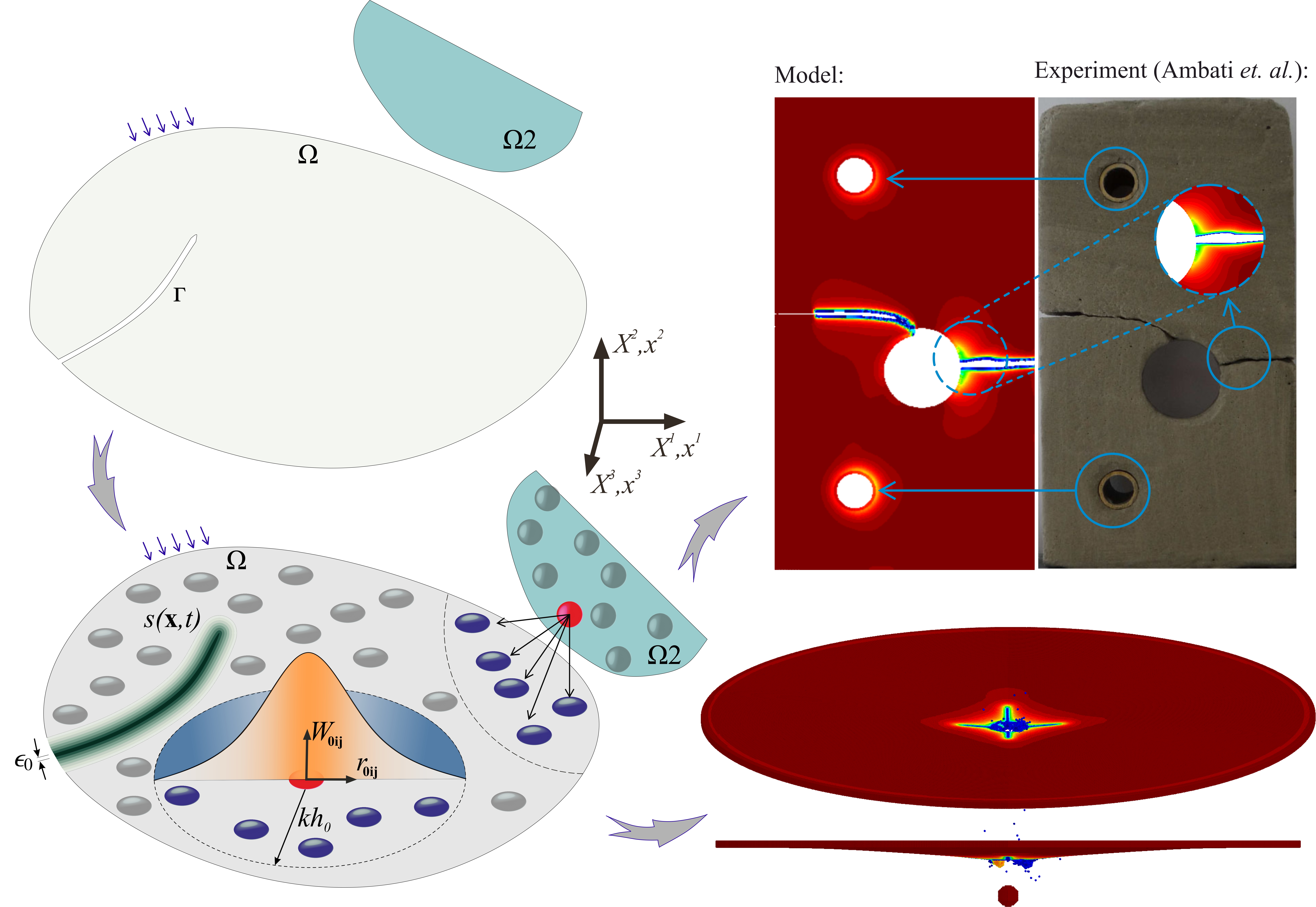}
\end{graphicalabstract}

\begin{highlights}
\item A framework for modeling brittle fracture within SPH is proposed.
\item SPH is posed in a total Lagrangian formulation to eliminate tensile instability.
\item Cracks are regularized over a length scale using a hyperbolic phase field approach.
\item The hyperbolic phase field equation is compatible with explicit time integration.
\item Challenging problems with complex fracture paths are presented.
\end{highlights}

\end{frontmatter}

\section{Introduction}
\label{sec:introduction}
The Smoothed Particle Hydrodynamics (SPH) is one of the oldest and well established particle methods, and has an outstanding application history in many fields. SPH was originally developed as an interpolation technique to study astrophysics related problems \cite{gingold1977smoothed,Lucy1977}, and was later extended to discretize fluid mechanics equations \cite{MONAGHAN1994399}. One of the major advantages of SPH over other existing particle methods, is that the inclusion of new physics is quite straightforward, and it can therefore be readily extended to new areas of application \cite{Libersky2006}. When it comes to solid mechanics, its ability to handle extreme material distortion makes it ideal for simulating large deformation problems, where conventional Lagrangian mesh-based techniques fail due to
mesh entanglement and the need for frequent mesh updating or remeshing. Many interesting
variants of the technique have been developed since it’s initial appearance, and the interested reader
should consult \cite{gingold1977smoothed,Lucy1977,MONAGHAN1994399,Swelge1995,Dyka1995,Bonet1999,Monaghan2005a,Libersky2006,Reveles2007,liu2010smoothed,Sigalotti2016,Chen2016,Rahimi2022114370} and the references therein.

Although SPH has been successfully applied to some complex problems in solid mechanics, fracture is still a challenging and open research area within the SPH community, and has not been extensively studied. This is because SPH is traditionally used to discretize equations of continuum mechanics, thus modeling of fracture and distinct crack surfaces is difficult. SPH researchers have so far relied on ad-hoc empirical local damage models, ``cracking particles" approaches, pseudo-spring and virtual link approaches, and cohesive zone models. All the aforementioned techniques have undoubtedly advanced the state-of-the-art in simulating fracture within SPH. However, they come with certain limitations and drawbacks. For example, local damage \cite{rabczuk2003simulation} models are based on empirical damage laws rather than comprehensive fracture theories, and lead to mesh dependency and non-convergent results under refinement. The ``cracking particles" approach \cite{Rabczuk2004a} resembles the extended finite element method (XFEM) \cite{moes1999finite,moes2002extended}, therefore the fracture surfaces need to be tracked, whereas the local enrichment of the approximation space leads to increased computational costs. In the pseudo-spring approach \cite{chakraborty2013pseudo,islam2019total,islam2020pseudo} the damage evolution is based on rather simple linear damage models, the softening curve of the damage law might lead to instabilities, and past research has shown that they are prone to spurious damage patterns. Finally, when it comes to cohesive zone models, the authors in \cite{wang2019new,wang2020simulation,bui2021smoothed} developed a continuum constitutive model featuring a cohesive fracture process zone, with an intrinsic physical length scale that helps remove the mesh dependency, while at the same time fracture surfaces do not need to be explicitly tracked. Nonetheless, kinematic enrichment is performed (similar to XFEM), which might potentially lead to increased computational times and complex implementation. Even though the aforementioned crack simulation techniques proved sufficient for certain classes of problems, it is evident that there is significant room for improvement when it comes to modeling fracture within the SPH framework.

Recently, phase field models of brittle fracture gained popularity within academic circles. Many interesting variants of the method have emerged and several challenging applications have been addressed \cite{kuhn2008phase,miehe2010phase,Borden2012,ambati2015review,wu2018length,Egger2019,wu2020phase,svolos2020thermal,mandal2021fracture,svolos2021anisotropic}. In phase field for fracture, cracks are not explicitly introduced in the solid, but instead the fracture surface is approximated by a phase field parameter that diffuses the discontinuity over a small region. The phase field parameter represents the material integrity, and it is a continuous variable that describes the smooth transition from the fully intact to the fully damaged state. The evolution of phase field, and hence of the fracture surface, is governed by a partial differential equation (PDE). Thus, a degree of non-locality is introduced, since the damage state at one point depends (through derivative information) on the states of neighboring points. The non-locality of the method leads to a well-posed mathematical model that exhibits mesh independence and convergence under refinement. At the same time, it can easily handle complicated discontinuity scenarios, such as crack initiation, propagation, coalescence, and branching, without prior knowledge of the crack paths, and without the need to embed any evolving discontinuities in the displacement field. As a result, propagating cracks are tracked automatically through the solution of the phase field PDE, which makes the method very attractive over other numerical approaches that require the explicit or implicit tracking of the discontinuities. Several approaches to phase field modeling of brittle fracture have been independently developed within the physics and mechanics communities, and for a thorough review the reader is encouraged to consult \cite{ambati2015review}. In the mechanics community, the starting point for deriving the method is the variational formulation of brittle fracture \cite{francfort1998revisiting}, which was regularized in \cite{bourdin2000numerical}, and extended Griffith’s theory of fracture \cite{griffith1921vi}.

Despite its many advantages, phase field modeling of brittle fracture has been mostly applied within the context of mesh-based numerical techniques, such as the finite element method (FEM) and isogeometric analysis (IGA) \cite{kuhn2008phase,miehe2010phase,Borden2012,borden2014higher,ambati2015review,chen2020adaptive,nguyen2021hybrid}. Employing particle and meshfree methods to discretize phase field equations for fracture comes with its own set of challenges, and very few relevant works can be found in the literature \cite{amiri2016fourth,kakouris2017phase,moutsanidis2018hyperbolic,kakouris2019phase,nguyen2019mesh,shao2019adaptive,li2020phase,wu2020efficient,shao2021adaptive}. To the best of the authors' knowledge, phase field models of brittle fracture have not been developed and investigated within SPH. In this proof-of-concept paper we develop a phase field model of brittle fracture suitable for SPH, and we demonstrate that phase field is a viable option to model fracture without running into issues that many of the previously used methods exhibit.

This paper is outlined as follows. In Section \ref{sec:SPH_TLSPH}, we review the basics of SPH, with a particular emphasis on total Lagrangian SPH. In Section \ref{sec:PhaseField} we provide the basics of phase field modeling of brittle fracture, with a focus on hyperbolic phase field models that are amenable to explicit time integration. In Section \ref{sec:TLSPH_Phasefield} we outline the proposed SPH framework for coupling phase field for fracture with solid mechanics. Section \ref{sec:results} contains several challenging numerical examples. Section \ref{sec:conclusion} draws conclusions and outlines future research directions.

\section{SPH approximation and the Total Lagrangian SPH}
\label{sec:SPH_TLSPH}
\subsection{Conventional SPH approximation}
\label{sec:ConvSPHapproach}
The SPH approach was initially proposed as a smoothed interpolation technique to deal with problems related to astrophysics \cite{gingold1977smoothed,Lucy1977}. The fundamental principle of SPH requires the integral representation of functions. The SPH utilizes a weighted interpolation technique to approximate the value of any arbitrary function $f \equiv f(\textbf{x})$ at any arbitrary point $\textbf{x}$ in Euclidean space. This is also known as the kernel approximation and can be written as \cite{gingold1977smoothed,Lucy1977}
\begin{equation}
	\label{eq:SPH_basic_interpolation}
	f(\textbf{x})\approxeq\int\displaylimits_{\Omega_{\textbf{x}}} f(\textbf{x}')~W(r,h)~d\textbf{x}' ,
\end{equation}
where $\textbf{x}'\in \Omega_{\textbf{x}}$ refers to the spatial coordinates of all the points located within the interpolation space $\Omega_{\textbf{x}}$ of point $\textbf{x}$. The term $W\equiv W(r,h)$, known as the kernel, represents the interpolation weights as a function of the Euclidean distance $r=|\textbf{x}-\textbf{x}'|$ and the smoothing length $h=1.33\Delta x$, where $\Delta x$ is the initial particle spacing. The interpolant in \prettyref{eq:SPH_basic_interpolation} reproduces the function $f$ exactly if the kernel is a delta function; that is $W=1$ if $\textbf{x}'=\textbf{x}$, and $W=0$ if $\textbf{x}'\neq \textbf{x}$. In practice, the kernel is chosen to be a compactly supported function and approaches the delta function as $h\rightarrow 0$.

The essential advantage of the SPH interpolation is that it allows an exact differentiation of the interpolant to produce the derivatives of function $f$. The spatial derivatives of $f$ can be computed using the exact differentiation of the kernel as \cite{gingold1977smoothed,Lucy1977}

\begin{equation}
	\label{eq:SPH_derivative}
	\nabla f\approxeq\int\displaylimits_{\Omega_{\textbf{x}}} f(\textbf{x}')~\nabla W~d\textbf{x}' \text{.}
\end{equation}
However, this kind of approximation does not vanish if the function is constant. To improve the accuracy of SPH for constant fields, \prettyref{eq:SPH_derivative} is rewritten by replacing $f$ with $\Phi f$, allowing for the nabla operator to be written as $\nabla (\Phi f) = \Phi \nabla f + f \nabla \Phi $. Thus, the derivative of the function becomes

\begin{equation}
    \label{eq:nabla_op}
	\nabla f = \frac{1}{\Phi} \left [\nabla (\Phi f) - f \nabla \Phi\right ] ,
\end{equation}
where $\Phi$ is any differentiable field. Generally, in continuum mechanics related applications, $\Phi$ is set to be the density $\rho$ since it represents a physical quantity in the continuity equation. By doing so, one arrives at the following form of SPH approximation which delivers a zero derivative for constant functions \cite{Monaghan2005a}

\begin{equation}
	\label{eq:derivative_monaghan}
	\nabla f(\textbf{x})\approxeq \frac{1}{\rho(\textbf{x})}\int\displaylimits_{\Omega_{\textbf{x}}} \rho(\textbf{x}') \left [ f(\textbf{x}') - f(\textbf{x}) \right ] ~\nabla W~d\textbf{x}'\text{.}
\end{equation}
In a spatially discretized SPH particle domain \prettyref{eq:SPH_basic_interpolation} and \prettyref{eq:derivative_monaghan} are written as

\begin{equation}
	\label{eq:SPH_disc}
	f_{\textbf{i}}\approxeq
	\sum_{\textbf{j}=1}^{N_{\textbf{i}}} f_{\textbf{j}}~W_{\textbf{ij}}~V_{\textbf{j}},
\end{equation}

\begin{equation}
	\label{eq:SPH_der_disc} 
	\frac{\partial f_{\textbf{i}}}{\partial x_{\textbf{i}}^{s}}
	\approxeq \frac{1}{\rho_{\textbf{i}}}
	\sum_{\textbf{j}=1}^{N_{\textbf{i}}} \rho_{\textbf{j}}~ (f_{\textbf{j}}-f_{\textbf{i}})~\frac{\partial W_{\textbf{ij}}}{\partial x_{\textbf{j}}^s} ~V_{\textbf{j}},
\end{equation}
in which $N_{\textbf{i}}$ is the total number of particles located within the interpolation space, also known as the neighborhood, influence, or support domain, of particle $\textbf{i}$. $V_{\textbf{j}}$ is the infinitesimal volume of particle $\textbf{j}$. $W_{\textbf{ij}}$ is the kernel function relating the particles $\textbf{i}$ and $\textbf{j}$, $x_{\textbf{j}}^s$ is the $s$ component of the spatial (or Eulerian) coordinate of particle $\textbf{j}$, and $\rho_{\mathbf{j}}$ is the density of particle $\textbf{j}$.

\begin{figure}[!htbp]
	\begin{center}
		\includegraphics[width=0.6\linewidth]{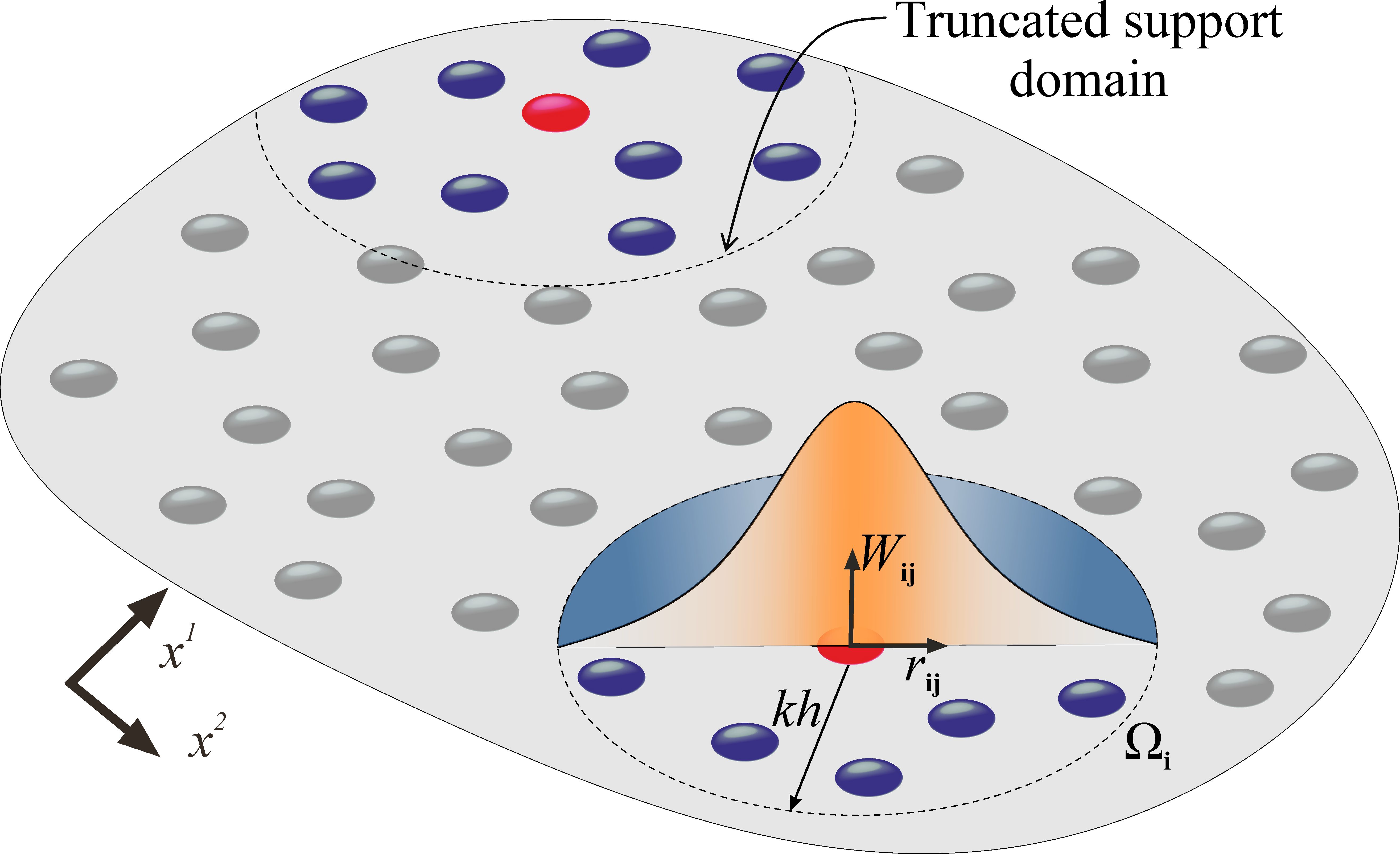}
	\end{center}
	\caption{\centering Discretized SPH particle domain with the illustration of full and truncated supports.}
	\label{fig:SPH_domain}
\end{figure}

As noted in \cite{Krongauz1997,Bonet1999,Sigalotti2016},
the discretization scheme introduced in \prettyref{eq:SPH_disc} and \prettyref{eq:SPH_der_disc} suffers from the particle inconsistency arising from the non-homogeneous distribution of particles and truncation of the support domain near the boundaries, as seen in \prettyref{fig:SPH_domain}. This leads to a lack of conservation of mass, and linear and angular momentum. 
A number of correction strategies to restore the particle consistency based on the kernel and gradient correction have been introduced \cite{Bonet1999,Chen2016,Vila1999}. For its simplicity, throughout this paper, we follow the gradient correction approach presented in \cite{Bonet1999}, in which the gradient is corrected as $\Tilde{\nabla}W_{\textbf{ij}}=\textbf{C}_{\textbf{i}}\cdot \nabla W_{\textbf{ij}}$.
\begin{equation}
	\label{eq:Correction_tens} 
	\textbf{C}_{\textbf{i}}=-\left(   \sum_{\textbf{j}=1}^{N_{\textbf{i}}}  \left( \nabla W_{\textbf{ij}} \otimes \textbf{r}_{\textbf{ij}}\right) V_{\textbf{j}}      \right)^{-1}\text{,}
\end{equation}
is a second rank correction tensor, and $\textbf{r}_{\textbf{ij}}=\textbf{x}_{\textbf{i}} - \textbf{x}_{\textbf{j}}$ is the relative position vector of particles $\textbf{i}$ and $\textbf{j}$. The use of $\textbf{C}_{\textbf{i}}$ ensures that the gradient of any linear velocity field is exactly evaluated. From this point forward, for brevity, we will omit the approximation sign and will use $``\nabla W_{\textbf{ij}}"$ or $``\partial W/\partial \textbf{x}"$ to refer to the corrected form of the kernel gradient $\Tilde{\nabla}W_{\textbf{ij}}=\textbf{C}_{\textbf{i}}\cdot \nabla W_{\textbf{ij}}$.

\subsection{SPH in Total Lagrangian form}
So far, the general procedure for SPH as an approximation technique was presented. In this section, we further explore its applicability for solid mechanics related applications. It is well established that, when SPH is applied to solve the governing equations of solid mechanics, the use of Eulerian kernels (kernels defined on Eulerian coordinates) leads to the so-called tensile instability \cite{Swelge1995,Belytschko2000}. Such kind of instability is observed in the form of material distortion under tensile stress state. This phenomenon was first noticed in \cite{Swelge1995} by carrying out the Von Neumann stability analysis (see \cite{Swelge1995} for detailed explanation). The authors in \cite{Swelge1995} concluded that for a stable solution the following condition should be satisfied 
\begin{equation}
	\label{eq:stab_cond} 
	\frac{\partial^2 W}{\partial \textbf{x}^2} \pmb{\sigma} \leq 0,
\end{equation}
where $\pmb{\sigma}$ is the stress. This implies that as long as the left hand side of \prettyref{eq:stab_cond} is greater than 0 for a particle, the SPH approximation delivers unstable solutions. Several remedies have been suggested to overcome such instability. For example, \cite{Morris1995} suggested an approach based on ``Mutating the kernel" (see Section 6.6 of \cite{Morris1995}). Since the tensile instability has a close relation with the second derivative of the kernel, using a proper kernel will result in a stable solution. However, such an approach is only effective in special cases \cite{Reveles2007}. Instead, the authors in \cite{Dyka1995} suggested the use of ``stress points" in which additional computational nodes are introduced away from the original SPH particles to carry the stress information separately. However, in \cite{Rabczuk2004} it was noted that the use of Eulerian kernels along with the stress points does not fully mitigate the innate instability of SPH. The authors in \cite{Rabczuk2004} proposed the use of Lagrangian kernels which deliver a more stable solution. The difference between Eulerian and Lagrangian kernels is that the former is a function of Eulerian (or spatial) coordinates, whereas the later is established based on the Lagrangian/reference (or material) coordinates.

\begin{figure}
	\begin{center}
		\includegraphics[width=\linewidth]{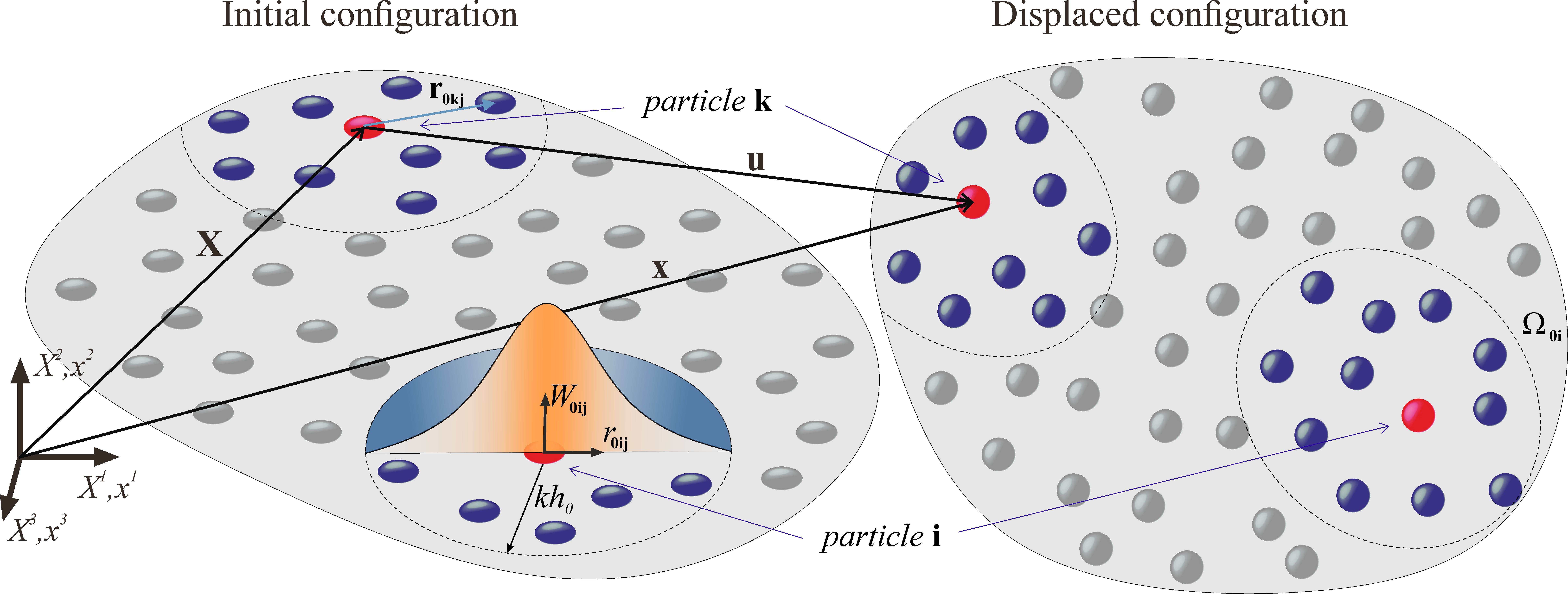}
	\end{center}
	\caption{\centering Initial and displaced configurations in TLSPH.}
	\label{fig:TLSPH_domain}
\end{figure}

In the Total Lagrangian SPH (TLSPH) formalism, the kernel function and its derivatives are written as functions of the Lagrangian coordinates, thus, leading to the TLSPH approximation of an arbitrary function and its derivative as

\begin{equation}
	\label{eq:TLSPH_disc}
	f_{\textbf{i}}=
	\sum_{\textbf{j}=1}^{N_{\textbf{i}}} \frac{m_{0\textbf{j}}}{\rho_{0\textbf{j}}}~f_{\textbf{j}}~W_{\textbf{0ij}},
\end{equation}
\begin{equation}
	\label{eq:TLSPH_der_disc} 
	\frac{\partial f_{\textbf{i}}}{\partial X_{\textbf{i}}^{s}}
	= \frac{1}{\rho_{0\textbf{i}}}
	\sum_{\textbf{j}=1}^{N_{\textbf{i}}} m_{0\textbf{j}}~ (f_{\textbf{j}}-f_{\textbf{i}})~\frac{\partial W_{0\textbf{ij}}}{\partial X_{\textbf{j}}^s},
\end{equation}
where
\begin{equation}
	\label{eq:disp}
	\textbf{X}=\textbf{x}-\textbf{u}.
\end{equation}
$W_{0\textbf{ij}}$, $m_0$, and $\rho_0$ are the kernel, mass, and density of the associated particle evaluated in the reference (initial) configuration, respectively. $\textbf{u}$ is the displacement, and $\textbf{X}$ is the Lagrangian coordinate, as depicted in \prettyref{fig:TLSPH_domain}. 
One major advantage of TLSPH is that, unlike the conventional SPH, its support domain does not change and is evaluated only once at the beginning of the simulation. This decreases the computational cost considerably. However, it becomes problematic in situations where the domain distortion is such that it leads to significant changes in a particle's neighbors. Cases like this are not considered in this paper and will be studied in subsequent works.

\subsection{Governing Equations in TLSPH framework}
\label{sec:GovEqTLSPH}

In this section we present the governing equation of elastic dynamics in TLSPH formalism. Our primary concern is the conservation of momentum, given that the mass is constant, and the total energy is naturally conserved\footnote{Strictly speaking, due to the presence of artificial viscosity, the energy is not exactly conserved. However, the small difference is assumed to be negligible for practical purposes. For more details on the aspect of artificial viscosity the interested reader should consult \cite{MONAGHAN198571}.}. In elastodynamics, the momentum balance equation in the reference configuration is written as \cite{Rabczuk2004}
\begin{equation}
	\label{eq:momentumEquation}
	\frac{d \textbf{v}}{dt} = \frac{1}{\rho_{0}}\nabla_{0} \cdot \textbf{P} + \textbf{b}_{0} \; \text{in} \; \Omega_{0},
\end{equation}
where $\Omega_{0}$ is the domain of the continuum body in the reference coordinates, and $\textbf{v}$, $\textbf{b}_0$, and $t$, are the velocity, external force, and time, respectively. $\nabla_{0} \cdot \textbf{P}$ is the divergence of the first Piola--Kirchhoff stress tensor with respect to the reference coordinates. The first Piola--Kirchhoff stress tensor, $\textbf{P}$, is computed through an appropriate constitutive model, in which the deformation gradient is the corresponding measure of deformation. Details of the specific constitutive models employed in this work will be presented in subsequent sections. The deformation gradient is given as

\begin{equation}
	\label{eq:deformGrad}
	\textbf{F} = \frac{\partial\textbf{x}}{\partial\textbf{X}}=\frac{\partial(\textbf{u}+\textbf{X})}{\partial\textbf{X}} = 
	\frac{\partial\textbf{u}}{\partial\textbf{X}}+\textbf{I},
\end{equation}
where $\textbf{I}$ is the identity matrix. In order to calculate the deformation gradient for particle $\textbf{i}$ we utilize the TLSPH approach as follows
\begin{equation}
	\label{eq:deformGradTLSPH}
	\textbf{F}_\textbf{i} = \textbf{I}+ \frac{1}{\rho_{0\textbf{i}}}
	\sum_{\textbf{j}=1}^{N_{\textbf{i}}} m_{0\textbf{j}}~ (\textbf{u}_{\textbf{j}}-\textbf{u}_{\textbf{i}})\otimes \nabla_0 W_{0\textbf{ij}},
\end{equation}
or in index notation as
\begin{equation}
	\label{eq:deformGradTLSPHindex}
	F^{ks}_\textbf{i} = \delta^{ks} + \frac{1}{\rho_{0\textbf{i}}}
	\sum_{\textbf{j}=1}^{N_{\textbf{i}}} m_{0\textbf{j}}~ u^{k}_{\textbf{ji}}~\frac{\partial W_{0\textbf{ij}}}{\partial X_{\textbf{j}}^s},
\end{equation}
in which $u_{\textbf{ji}}^k$ is the $k$ component of the displacement difference vector $\textbf{u}_{\textbf{ji}}=\textbf{u}_{\textbf{j}}-\textbf{u}_{\textbf{i}}$, and $\delta^{ks}$ is the Kronecker delta. Recall that $\partial W_{0\textbf{ij}}/\partial X_{\textbf{j}}^s$ is the corrected derivative of the Lagrangian kernel (see \prettyref{sec:ConvSPHapproach}). Similarly, \prettyref{eq:momentumEquation} can be written for particle $\textbf{i}$ as \cite{Lin2015}
\begin{equation}
	\label{eq:MomentumEqTLSPH}
	\frac{d \textbf{v}_\textbf{i}}{dt} = \sum_{\textbf{j}=1}^{N_{\textbf{i}}} 
	m_{0\textbf{j}}
	\left( 
	\frac{\textbf{P}_\textbf{i}}{\rho^2_{0\textbf{i}}} +
	\frac{\textbf{P}_\textbf{j}}{\rho^2_{0\textbf{j}}} +
	\textbf{P}_{v\textbf{ij}}
	\right)
	\cdot \nabla_0 W_{0\textbf{ij}} + \textbf{b}_{0\textbf{i}},
\end{equation}
or in index notation as
\begin{equation}
	\label{eq:MomentumEqTLSPHindex}
	\frac{d v^k_\textbf{i}}{dt} = \sum_{\textbf{j}=1}^{N_{\textbf{i}}} 
	m_{0\textbf{j}}
	\left( 
	\frac{P^{ks}_\textbf{i}}{\rho^2_{0\textbf{i}}} +
	\frac{P^{ks}_\textbf{j}}{\rho^2_{0\textbf{j}}} +
	P^{ks}_{v\textbf{ij}}
	\right)
	~\frac{\partial W_{0\textbf{ij}}}{\partial X_{\textbf{j}}^s} + b_{0\textbf{i}}^k,
\end{equation}
where Einstein's summation rule is employed for the repeated index $s$. Furthermore, the artificial viscosity term $\textbf{P}_{v\textbf{ij}} = \textrm{det}(\textbf{F}_{\textbf{i}})~\pi_{\textbf{ij}}~\textbf{F}_{\textbf{i}}^{-1}$ is included to avoid numerical instabilities arising from zero-energy mode discrepancy in the form of a jump in the field variables or shock waves \cite{He2017,islam2019total,monaghan1983374,Rahimi2022114370}. Following the work of \cite{monaghan1983374} we compute the coefficient $\pi_{\textbf{ij}}$ as
\begin{equation}
    \label{eq:artViscCoef}
	\pi_{\textbf{ij}}=\frac{1}{\rho_{0\textbf{i}}}(\beta_2 G_{\textbf{ij}}^2-\beta_1~c_{0\textbf{i}}G_{\textbf{ij}}),
\end{equation}
in which 
\begin{equation}
    \label{eq:sppedofsound}
	c_{0\textbf{i}}=c_{0}=\sqrt{\frac{\kappa+\frac{4}{3}\mu}{\rho_{0}}}
\end{equation}
is the speed of sound at particle $\textbf{i}$,
$\kappa = E/(3(1-2\nu))$ is the bulk modulus, $\mu=E/(2(\nu+1))$ is the shear modulus, with $E$ and $\nu$ being the elastic modulus and Poisson's ratio, respectively. $\beta_1$ and $\beta_2$ are scalar factors, and $G_{\textbf{ij}}$ is calculated as
\begin{equation}
	G_{\textbf{ij}}=\frac{h~(\textbf{v}_{\textbf{i}}-\textbf{v}_{\textbf{j}}) \cdot (\textbf{X}_{\textbf{i}}-\textbf{X}_{\textbf{j}})}{r^2_{0\textbf{ij}}+0.001 h^2}.
\end{equation}
Here, $r_{0\textbf{ij}}=|\textbf{X}_{\textbf{i}}-\textbf{X}_{\textbf{j}}|$ is the initial distance between particles $\textbf{i}$ and $\textbf{j}$. In most of the numerical examples in this paper we set $\beta_1 = 0.04$ and $\beta_2 = 0$. Note that for larger values of time step and lower particle resolutions, higher $\beta_1$ and $\beta_2$ values may be necessary. It should be mentioned that the artificial viscosity introduces nonphysical dumping to the system. Thus, the smallest possible values of $\beta_1$ and $\beta_2$ are preferred.

\section{Griffith's theory of brittle fracture and the phase field approximation}
\label{sec:PhaseField}
\begin{figure}[!htbp]
\begin{center}
	\begin{subfigure}{0.49\linewidth}
		
        \includegraphics[width=1\linewidth]{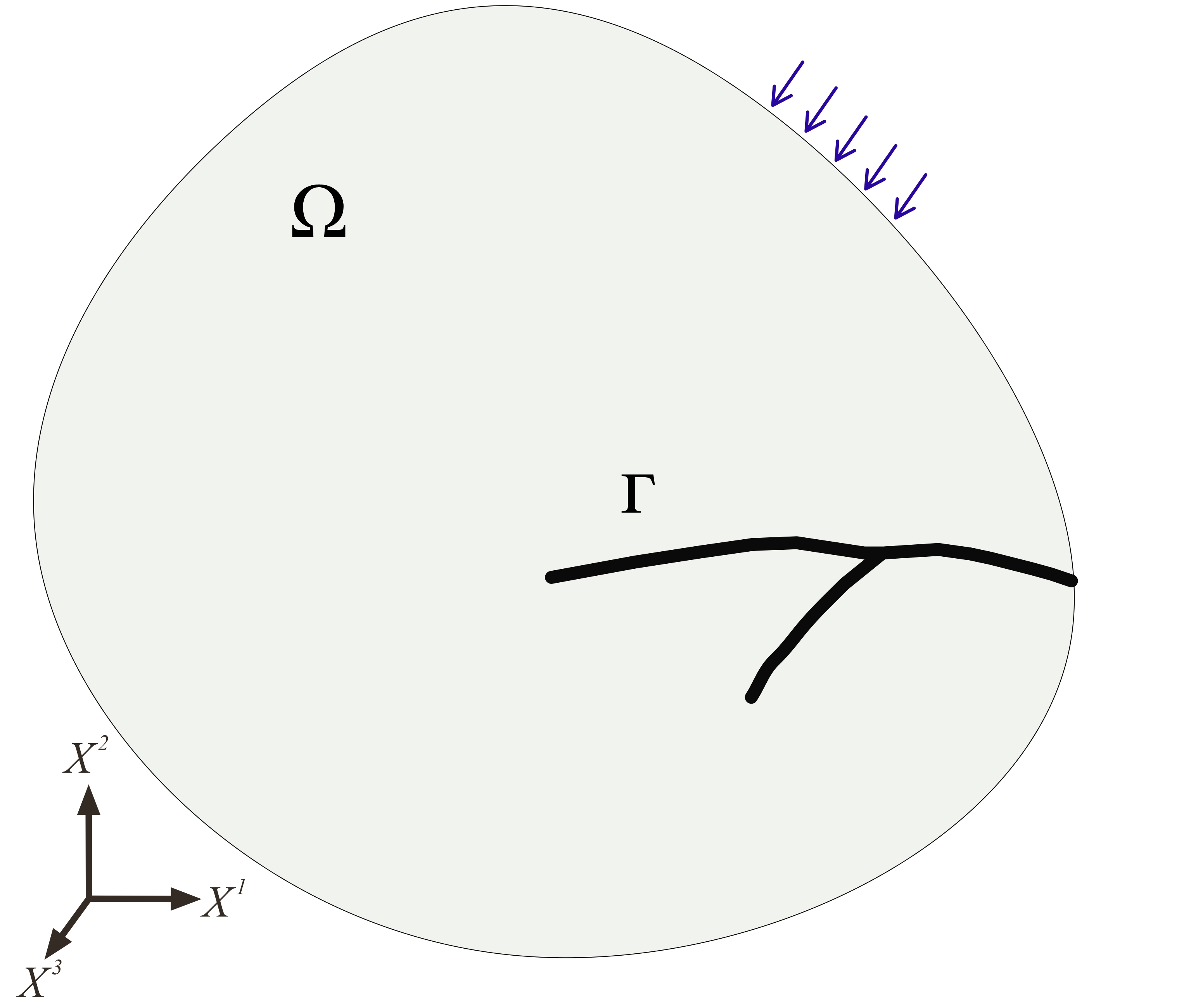}
 		\caption{}
	\end{subfigure}
	\begin{subfigure}{0.49\linewidth}
        \includegraphics[width=1\linewidth]{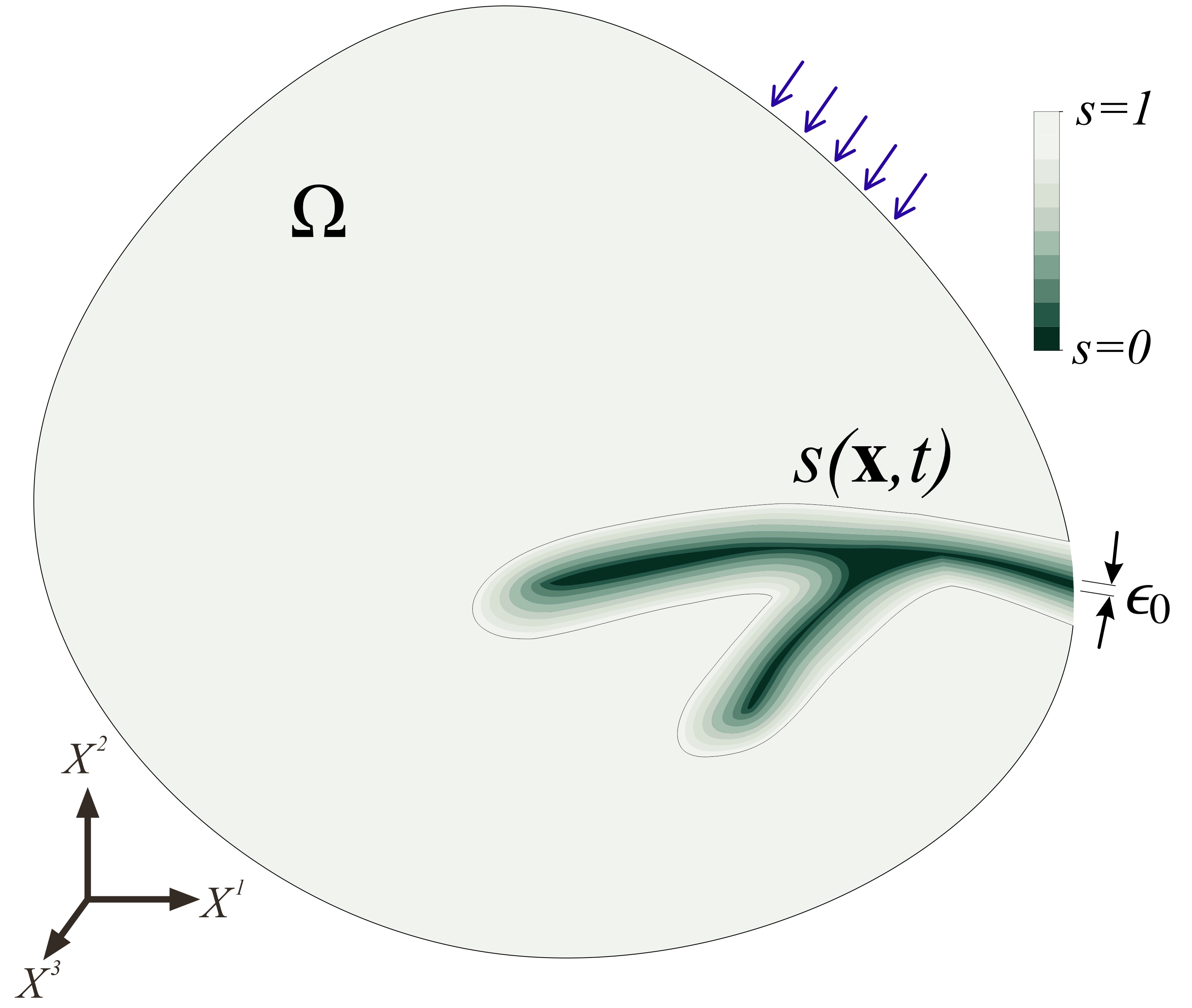}
		
 		\caption{}
	\end{subfigure}
	\caption{(a) A solid body $\Omega$ with internal discontinuity $\Gamma$. (b) Approximation of the internal discontinuity by the phase field $s(\mathbf{x},t)$. $\epsilon_{0}$ is the length scale parameter that controls the width of the discontinuity.}
	\label{fig:phasefield_domain}
	\end{center}
\end{figure}

In this section, to be self contained, we briefly present the basics of phase field modeling of brittle fracture. In the mechanics community, the prevailing models mostly originate from the variational formulation of brittle fracture \cite{francfort1998revisiting} and the related regularized formulation \cite{bourdin2000numerical,ambati2015review}. According to Griffith’s theory of brittle
fracture \cite{griffith1921vi}, the fracture energy density $G_{c}$ is the amount of energy required to open a unit area of crack surface. Then, the total potential energy of an elastic body $\Omega$, being the sum of the elastic energy
and the fracture energy, is given by the expression
\begin{equation}
\label{eq:energy_functional}
    \Pi_{p}(\mathbf{u},\Gamma) = \int_{\Omega} \psi_{e} (\bm{\varepsilon}({\mathbf{u}})) d \mathbf{x} + \int_{\Gamma} G_{c} d \Gamma,
\end{equation}
where $\psi_{e}$ is the elastic strain energy density, the details of which will be presented in subsequent sections, $\pmb{\varepsilon}$ is the infinitesimal strain tensor, and $\Gamma$ is the evolving internal discontinuity boundary which represents a set of discrete cracks. Extending this approach to dynamic problems the kinetic energy is defined as
\begin{equation}
\Pi_{k}(\mathbf{\dot{u}}) = \frac{1}{2} \int_{\Omega} \rho~\dot{\textbf{u}}\cdot \dot{\textbf{u}}~d \mathbf{x},
\end{equation}
where $\rho$ is the material density in the current configuration, and the superimposed dot denotes time differentiation. Combining the kinetic energy with the potential energy of \prettyref{eq:energy_functional} we arrive at the Lagrangian of the discrete fracture problem
\begin{equation}
\label{eq:Lagrangian}
L (\mathbf{u}, \mathbf{\dot{u}}, \Gamma) = \Pi_{k}(\mathbf{\dot{u}}) - \Pi_{p}(\mathbf{u},\Gamma) = \int_{\Omega} \Bigg ( \frac{1}{2} \rho ~\dot{\textbf{u}}\cdot \dot{\textbf{u}}~ - \psi_{e} (\pmb{\varepsilon} (\mathbf{u})) \Bigg ) d \mathbf{x} - \int_{\Gamma} G_{c} d\Gamma. 
\end{equation}
The Euler-Lagrange equations of this functional determine the equations of motion of the body, including the entire process of crack initiation, propagation, and branching of preexisting cracks. However, the numerical treatment is quite complex because tracking of the evolving discontinuity $\Gamma$ is required, which leads to complicated and expensive computations \cite{kuhn2008phase,Borden2012}. Therefore, to circumvent the above-mentioned difficulties, the regularized expression for the fracture energy was instead proposed in \cite{bourdin2000numerical} given as
\begin{equation}
\label{eq:fracture_regularized}
\int_{\Gamma} G_{c} d\Gamma \approx \int_{\Omega} G_{c} \Bigg ( \frac{(1-s)^2}{4 \epsilon_0} + \epsilon_0 | \nabla s |^2 \Bigg ) \, d\mathbf{x},
\end{equation}
where $s$ is the so-called phase field parameter (or damage parameter) which represents the fracture surface $\Gamma$. It is a continuous variable that describes the smooth transition from the undamaged state ($s = 1$) to the fully damaged one ($s = 0$). $\epsilon_0$ is a parameter that has dimension of a length and controls the width of the
smooth approximation of the crack (see \prettyref{fig:phasefield_domain}). When $\epsilon_0 \to 0$ the phase field approximation converges to the discrete fracture surface. In practice, $\epsilon_0$ should be sufficiently small so that the physics of the problem are not violated, but at the same time it should be greater or roughly equal to the spatial discretization size so that the crack is regularized over some finite width. To model the loss of material stiffness, the elastic strain energy density is defined as \cite{miehe2010phase,Borden2012}
\begin{equation}
\label{eq:energy_density}
\psi_{e} (\pmb{\varepsilon},s) = s^2 \psi_{e}^{+} + \psi_{e}^{-}.
\end{equation}
$\psi_{e}^{+}$ and $\psi_{e}^{-}$ are the positive and negative parts of the elastic strain energy density, respectively, which will be defined in the next section according to the particular constitutive models employed. As can become evident from \prettyref{eq:energy_density}, crack propagation is only allowed in tension since the phase field parameter is applied only to the tensile part of the elastic strain energy. By substituting the phase field approximations for the fracture energy \prettyref{eq:fracture_regularized} and the elastic energy density \prettyref{eq:energy_density} into the expression for the Lagrangian functional \prettyref{eq:Lagrangian}, and by deriving the corresponding Euler-Lagrange equations, we arrive at the strong form equations of motion that consist of the momentum balance presented in \prettyref{eq:momentumEquation} and the phase field parameter evolution equation
\begin{equation}
\label{eq:phasefield_evolution1}
\bigg ( \frac{4 \epsilon_0 \psi_{e}^{+}}{G_{c}} + 1 \bigg ) s - 4 \epsilon_0^{2} \nabla^2 s = 1 \; \text{in} \; \Omega_{0},
\end{equation}
where $\nabla^2 s$ is the Laplacian of $s$.
\begin{remark}
It can be seen that the kinetic energy term in the Lagrange energy functional is not affected by the phase field parameter $s$, which leads to conservation of mass.
\end{remark}
In order to model the irreversibility condition (cracks do not heal) a history functional is introduced as 
\begin{equation}
\label{eq:irreversibility}
\mathcal{H} (\mathbf{X},t) = \max_{\tau<=t} ( \psi_{e}^{+} (\mathbf{X},\tau)),
\end{equation}
and is used in place of the tensile elastic strain energy $\psi_{e}^{+}$ in the phase field's governing equation. The updated strong form that governs the evolution of the phase field parameter is
\begin{equation}
\label{eq:phasefield_evolution2}
\bigg ( \frac{4 \epsilon_0 \mathcal{H}}{G_{c}} + 1 \bigg ) s - 4 \epsilon_0^{2} \nabla^2 s = 1 \; \text{in} \; \Omega_{0}.
\end{equation}
Finally, initial conditions need to be introduced for the history functional, i.e.
\begin{equation}
\label{eq:phasefield_initial}
\mathcal{H}(\mathbf{X},0) = \mathcal{H}_{0} (\mathbf{X}).
\end{equation}
The initial condition of the history functional $\mathcal{H}_{0} (\mathbf{X})$ can be used to model preexisting cracks in the domain. \prettyref{eq:phasefield_evolution2} and the corresponding initial condition \prettyref{eq:phasefield_initial} need to be solved together with the equations of motion (e.g. momentum balance \prettyref{eq:momentumEquation}) to compute both the displacement field $\mathbf{u} (\mathbf{X},t)$ and the phase field parameter $s (\mathbf{X},t)$. 

Although the aforementioned approach and its variants have been successfully applied within the context of numerical methods based on implicit time integration, they become problematic when lumped-mass explicit dynamics schemes are employed. The reason is that \prettyref{eq:phasefield_evolution2} is an elliptic PDE, and therefore a linear system needs to be solved in every time step to determine the phase field parameter. Though it is possible from an implementation point of view, the solution of an elliptic equation in every explicit step would make the approach computationally expensive. To overcome this difficulty, phase field models which include time dependency have been proposed. In \cite{kuhn2008phase} the phase field parameter evolves according to a parabolic PDE as
\begin{equation}
\label{eq:phasefield_parabolic}
\frac{1}{M} \dot{s} + 2 s \psi_{e}^{+} - G_{c} \bigg ( 2 \epsilon_0 \nabla^2 s + \frac{1-s}{2 \epsilon_0} \bigg ) = 0 \; \text{in} \; \Omega_{0},
\end{equation}
where $M$ is a parameter controlling the rate at which local damage information diffuses into the bulk material. Clearly, when $M \to \infty$ the model approaches the standard elliptic models like the one of \prettyref{eq:phasefield_evolution1}. However, \prettyref{eq:phasefield_parabolic} is a heat-equation-like PDE, and is subjected to an unfavorable CFL condition when explicit time integration is employed \cite{Kamensky2018}
\begin{equation}
\label{eq:parabolic_stability}
\Delta t \leq \Delta x^2,
\end{equation}
where $\Delta t$ and $\Delta x$ are time and length scales associated with the discretization. Therefore, the authors in \cite{Kamensky2018} proposed a novel model in which the phase field parameter evolves according to a hyperbolic PDE by adding a second-order time derivative to the model of \prettyref{eq:phasefield_parabolic}. The resulting PDE is
\begin{equation}
\label{eq:phasefield_hyperbolic1}
\frac{2 G_{c} \epsilon_0}{c^2} \ddot{s} + \frac{1}{M} \dot{s} + 2 s \psi_{e}^{+} - G_{c} \bigg ( 2 \epsilon_0 \nabla^2 s + \frac{1-s}{2 \epsilon_0} \bigg ) = 0 \; \text{in} \; \Omega_{0}.
\end{equation}
$c$ is a speed limit on the propagation of the phase field parameter through the undamaged material and is taken to be equal to the sound speed given in \prettyref{eq:sppedofsound}. The model of \prettyref{eq:phasefield_hyperbolic1} is subjected to the hyperbolic stability condition of
\begin{equation}
\label{eq:hyperbolic_stability}
c \Delta t \leq \Delta x,
\end{equation}
which is much less restrictive than the parabolic one (\prettyref{eq:parabolic_stability}), especially in the limit of $\Delta x \to 0$. To avoid any wave-like behavior due to the hyperbolic nature of the governing equation, the authors in \cite{Kamensky2018} proposed an upper bound on $M$ such that the system is overdamped and the behavior of the phase field parameter is monotonic. The upper bound is given as
\begin{equation} \label{eq:Mparameter}
M \leq \frac{c}{2 \sqrt{4 G_{c} \epsilon_0 \psi_{e} + G_{c}^2}}.
\end{equation}
Finally, the irreversibility condition is enforced in a way similar to the elliptic problem (\prettyref{eq:irreversibility}) and the updated hyperbolic governing equation for the phase field damage parameter becomes
\begin{equation}
\label{eq:phasefield_hyperbolic2}
\frac{2 G_{c} \epsilon_0}{c^2} \ddot{s} + \frac{1}{M} \dot{s} + 2 s \mathcal{H} - G_{c} \bigg ( 2 \epsilon_0 \nabla^2 s + \frac{1-s}{2 \epsilon_0} \bigg ) = 0 \; \text{in} \; \Omega_{0},
\end{equation}
with the initial condition of \prettyref{eq:phasefield_initial}, whereas the upper bound on the parameter M becomes
\begin{equation} \label{eq:Mparameter1}
M \leq \frac{c}{2 \sqrt{4 G_{c} \epsilon_0 \mathcal{H} + G_{c}^2}}.
\end{equation}
Traditionally, explicit time integration is used in SPH. Although there have been some efforts to develop SPH variants based on implicit time integration (e.g. \cite{peer2018implicit}), the process is rather tedious and computationally expensive due to the large number of neighbors associated with each particle, and requires the solution of large linear systems. Thus, the hyperbolic version of the phase field PDE (\prettyref{eq:phasefield_hyperbolic2}) is employed in this work.

\section{Proposed framework}
\label{sec:TLSPH_Phasefield}
\subsection{Governing equations}

With the two previous sections at hand (\prettyref{sec:SPH_TLSPH} and \prettyref{sec:PhaseField}), we present a total Lagrangian SPH framework for the solution of the coupled elastodynamics--phase field problem expressed through the following strong form governing equations,
\begin{equation}
  \left\{\begin{array}{lr}	\frac{d \textbf{v}}{dt} = \frac{1}{\rho_{0}}\nabla_{0} \cdot \textbf{P} + \textbf{b}_{0}&\text{in} \; \Omega_{0} \\\frac{2 G_{c} \epsilon_0}{c^2} \ddot{s} + \frac{1}{M} \dot{s} + 2 s \mathcal{H} - G_{c} \bigg ( 2 \epsilon_0 \nabla^2 s + \frac{1-s}{2 \epsilon_0} \bigg ) = 0 & \text{in} \; \Omega_{0} \end{array}\right. .
\end{equation}
In the above-mentioned coupled problem, the first equation corresponds to the balance of linear momentum, that was briefly outlined in Section~\ref{sec:SPH_TLSPH}, whereas the second equation corresponds to the phase field (damage) parameter evolution. The solution of the phase field PDE is based upon the basics of TLSPH, and the overall solution procedure is outlined in this section. Before proceeding to the solution procedure we provide the details of the constitutive model.

\subsection{Constitutive modeling}

In Section \ref{sec:SPH_TLSPH} we presented the basics of TLSPH for discretizing solid mechanics equations in terms of the first Piola--Kirchhoff stress, whereas in Section \ref{sec:PhaseField} we outlined the basics of the phase field approach for brittle fracture in terms of a general elastic strain energy density $\psi_{e}$. Here, we provide the details of the elastic strain energy density and the corresponding constitutive equations based on two hyperelastic models; the isotropic Saint-Venant Kirchhoff and the neo-Hookean model. However, we would like to point out that other constitutive models can be used with the proposed framework.

\subsubsection{Saint-Venant Kirchhoff}
\label{sec:const_SVK}
The Saint-Venant Kirchhoff model can be easily derived by extending the linear-elastic framework presented in \cite{Borden2012} to large rotations. This is achieved by replacing the infinitesimal strain tensor with the Green--Lagrange strain tensor $\mathbf{E} = \frac{1}{2} (\mathbf{F}^T \mathbf{F} - \textbf{I})$. Then, the elastic strain energy density functional becomes
\begin{equation}
\psi_{e} = \frac{1}{2} \lambda ( \text{tr}  \mathbf{E})^2 + \mu \text{tr} (\mathbf{E}^2),
\end{equation}
where $\lambda=\kappa -\frac{2}{3}\mu$ and $\mu=E/(2(\nu+1))$ are the Lame parameters.
We then define
\begin{equation}\label{eq:elasEnergyP}
    \psi_{e}^{+} = \frac{1}{2} \lambda \{ \text{tr}  \mathbf{E} \}_{+}^2 + \mu \text{tr} (\mathbf{E}^{+} \mathbf{E}^{+}),
\end{equation} 
\begin{equation}\label{eq:elasEnergyN}
    \psi_{e}^{-} = \frac{1}{2} \lambda \{ \text{tr}  \mathbf{E} \}_{-}^2 + \mu \text{tr} (\mathbf{E}^{-} \mathbf{E}^{-}),
\end{equation}
where the following decomposition is employed
\begin{equation}\label{eq:GLstratinP}
\mathbf{E}^{+} = \mathbf{Q} \mathbf{\Lambda}^{+} \mathbf{Q}^{T},
\end{equation}
\begin{equation}\label{eq:GLstratinN}
\mathbf{E}^{-} = \mathbf{Q} \mathbf{\Lambda}^{-} \mathbf{Q}^{T},
\end{equation}
\begin{equation}
\mathbf{E} = \mathbf{Q} \mathbf{\Lambda} \mathbf{Q}^{T}.
\end{equation}
$\Lambda = \text{diag} (\lambda_{1}, \lambda_{2}, \lambda_{3} )$ has the eigenvalues of $\mathbf{E}$ on its diagonal, $\mathbf{Q}$ has the corresponding eigenvectors as its columns, $\Lambda^{\pm} = \text{diag} (\lambda_{1}^{\pm}, \lambda_{2}^{\pm}, \lambda_{3}^{\pm} )$, and $\{\,.\,\}_{\pm}$ selects the $\pm$ part of its argument, i.e.
\begin{equation}\label{eq:xpm_decomp}
    \{x\}_{\pm} = \left\{\begin{array}{lr}x & x\in\mathbb{R}^\pm\\ 0 &\text{otherwise}\end{array}\right.\text{ .}
\end{equation}
The second Piola--Kirchhoff stress can then be computed by differentiating the strain energy density $\psi_{e}$ with respect to the Green--Lagrange strain tensor $\mathbf{E}$,
\begin{equation}
\label{eq:2PK_SVK1}
\mathbf{S}^{\pm} = \frac{\partial \psi_{e}^{\pm}}{\partial \mathbf{E}} = \lambda~\{ \text{tr}  \mathbf{E} \}_{\pm} \textbf{I}+2\mu \textbf{E}^{\pm},
\end{equation}
and
\begin{equation}
\label{eq:2PK_SVK2}
\mathbf{S} = s^2 \mathbf{S}^{+} + \mathbf{S}^{-}.
\end{equation}
Finally, the first Piola–Kirchhoff stress can be computed as 
\begin{equation}
\label{eq:PK-2PK}
\mathbf{P} = \mathbf{F} \mathbf{S}.
\end{equation}

\subsubsection{Neo-Hookean}
\label{sec:const_NH}
The Saint-Venant Kirchhoff model is known to exhibit instabilities in the case of strong compression \cite{holzapfel2002nonlinear}. Thus, for the problems in this paper involving strong compression, we employ a Neo-Hookean material. For the Neo-Hookean model we follow the presentation in \cite{borden2016phase,moutsanidis2018hyperbolic}. The positive and negative parts of the elastic strain energy density are given as
\begin{equation}\label{eq:ksipNH}
\psi_{e}^{+} = \left\{\begin{array}{lr}	U(J) + \overline{\psi}_{e} (\overline{\mathbf{b}}) \, & J \geq 1 \\ \overline{\psi}_{e} (\overline{\mathbf{b}}) \, & J< 1, \end{array}\right.
\end{equation}
\begin{equation}\label{eq:ksinNH}
\psi_{e}^{-} = \left\{\begin{array}{lr}	0 \, & J \geq 1 \\ U(J) \, & J< 1, \end{array}\right.
\end{equation}
where
\begin{equation} \label{eq:UjNeoHok}
U(J) = \frac{1}{2} \kappa~\bigg( \frac{1}{2} (J^2 - 1) - \text{ln} J \bigg),
\end{equation}
\begin{equation} \label{eq:ksiNeoHok}
\overline{\psi}_{e} (\overline{\mathbf{b}}) = \frac{1}{2} \mu \Big( \text{tr} \overline{\mathbf{b}} - 3 \Big),
\end{equation} 
\begin{equation}\label{eq:jacobian}
J = \text{det} \mathbf{F},
\end{equation}
\begin{equation}\label{eq:bmatrix}
\mathbf{b} = \mathbf{F} \mathbf{F}^T,
\end{equation}
\begin{equation}\label{eq:bbarmatrix}
\overline{\mathbf{b}} = J\,^{-2/3} \mathbf{b}.
\end{equation}
The second Piola–-Kirchhoff stress can then be computed as
\begin{equation}
\mathbf{S} = 2 \frac{\partial \psi_{e}}{\partial \mathbf{b}}.
\end{equation}
For the given elastic strain energy density this results in
\begin{equation}
\label{eq:2PK_NH1}
\mathbf{S} = 2 \left\{\begin{array}{lr} s^2 \bigg( U'(J) \frac{\partial J}{\partial \mathbf{b}} + \frac{\partial \overline{\psi}_{e} (\overline{\mathbf{b})}}{\partial \mathbf{b}} \bigg) \, & J \geq 1 \\ U'(J) \frac{\partial J}{\partial \mathbf{b}} + s^2 \frac{\partial \overline{\psi}_{e} (\overline{\mathbf{b})}}{\partial \mathbf{b}}  \, & J< 1. 
\end{array}\right.
\end{equation}
The derivatives in the above expression are computed as 
\begin{equation}
U'(J) = \frac{1}{2} \kappa \Big( J - J^{-1} \Big),
\end{equation}
\begin{equation}
\frac{\partial J}{\partial \mathbf{b}} = \frac{\partial \sqrt{\text{det} \mathbf{b}}}{\partial \mathbf{b}} = \frac{1}{2} J \mathbf{b}^{-1},
\end{equation}
\begin{equation}
\frac{\partial \overline{\psi}_{e}}{\partial \mathbf{b}} = \frac{\partial \overline{\psi}_{e}}{\partial \overline{\mathbf{b}}} \frac{\partial \overline{\mathbf{b}}}{\partial \mathbf{b}} = \frac{J^{-2/3}}{2} \mu \bigg( \mathbf{I} - \frac{1}{3} \big( \text{tr} \mathbf{b} \big) \mathbf{b}^{-1} \bigg).
\end{equation}
Substituting the above equations into \prettyref{eq:2PK_NH1} we get
\begin{equation}
\label{eq:2PK_NH2}
    \mathbf{S} = \left\{\begin{array}{lr} s^2 \bigg( \frac{1}{2} \kappa \big( J^2 - 1 \big) \mathbf{b}^{-1} + J^{-2/3} \mu \Big( \mathbf{I} - \frac{1}{3} \big( \text{tr} \mathbf{b} \big) \mathbf{b}^{-1} \Big) \bigg)  \, & J \geq 1 \\ \frac{1}{2} \kappa \big( J^2 - 1 \big) \mathbf{b}^{-1} +s^2 \bigg( J^{-2/3} \mu \Big( \mathbf{I} - \frac{1}{3} \big( \text{tr} \mathbf{b} \big) \mathbf{b}^{-1} \Big) \bigg)  \, & J< 1. 
\end{array}\right.
\end{equation}
Finally, the first Piola--Kirchhoff stress is computed as in \prettyref{eq:PK-2PK}.

\subsection{Solution procedure and numerical implementation}
\label{sec:procedure}

We first define two common terms that will be used regularly throughout this paper, the soft particle and the phase field limit.
\begin{definition}
    Soft (damaged) particle refers to any SPH particle whose phase field parameter drops sufficiently so that its stress is fairly low, and thus the effect on neighboring particles is negligible. From a numerical implementation point of view, soft particles behave as rigid objects having inertia but no internal force.
\end{definition}
\begin{definition}
    The phase field limit refers to the threshold value of the phase field parameter below which the particle is assumed to be a soft particle.
\end{definition}
In order to compute the first Piola--Kirchhoff stress through the constitutive model presented in the previous subsection, we redefine the deformation gradient as
\begin{equation} \label{eq:defGrad_phase}
    \textbf{F}_{\textbf{i}} = \left\{\begin{array}{lr}\frac{\partial\textbf{u}_{\textbf{i}}}{\partial\textbf{X}_{\textbf{i}}}+\textbf{I} & s_{\textbf{i}}>s_l\\ \textbf{I} &\text{otherwise}\end{array}\right.\text~{.}
\end{equation}
Here $s_l$ is the phase field limit that determines the particle's stiffness state and $\partial\textbf{u}_{\textbf{i}}/\partial\textbf{X}_{\textbf{i}}$ is calculated utilizing the corrected form of the TLSPH approach (see \prettyref{sec:GovEqTLSPH}) as
\begin{equation}
	\label{eq:displacementGrad}
	\frac{\partial u^k_{\textbf{i}}}{\partial X^s_{\textbf{i}}} = \frac{1}{\rho_{0\textbf{i}}}
	\sum_{\textbf{j}=1}^{N} m_{0\textbf{j}}~ u^{k}_{\textbf{ji}}~\frac{\partial W_{0\textbf{ij}}}{\partial X_{\textbf{j}}^s},
\end{equation}
where $u_{\textbf{ji}}^k$ is the $k$ component of the displacement difference vector $\textbf{u}_{\textbf{ji}}=\textbf{u}_{\textbf{j}}-\textbf{u}_{\textbf{i}}$ between particles \textbf{i} and \textbf{j}.
The kernel function in \prettyref{eq:displacementGrad}, and in any other calculation involving SPH interpolation, is chosen to be a cubic spline \cite{Monaghan1992} defined in reference coordinates as
\begin{equation}
	\label{eq:TLSPHkernel}
	W_{0\textbf{ij}}=\left\{\begin{array}{ll}0 & q>2\\ \frac{1}{4}C~(2-q)^3 &1\leq q\leq 2\\
	C~(1-1.5q^2+0.75q^3) & 0\leq q \leq 1
	\end{array}\right.\text~{,}
\end{equation}
in which $q=|\textbf{X}_{\textbf{i}}-\textbf{X}_{\textbf{j}}|/h$ and C is a constant determined as
\begin{equation}
	C=\left\{\begin{array}{ll}2/(3h) & \text{for one dimensional space}\\ 10/(7\pi h^2) & \text{for two dimensional space}\\
	1/(\pi h^3) & \text{for three dimensional space}
	\end{array}\right.\text~{,}
\end{equation}
where $h$ is the smoothing length. We further rearrange \prettyref{eq:phasefield_hyperbolic2} as
\begin{equation} \label{eq:phasefieldddot}
    \ddot{s} =\frac{c^2}{2 G_{c} \epsilon_0}  \left[ G_{c} \bigg ( 2 \epsilon_0 \nabla^2 s + \frac{1-s}{2 \epsilon_0} \bigg )- \frac{1}{M} \dot{s} - 2 s \mathcal{H} \right]
\end{equation}
in order to have an explicit definition of the phase field inertia $\ddot{s}$, in terms of the phase field parameter $s$, its time derivative $\dot{s}$, and its Laplacian $\nabla^2 s$. The Laplacian of the phase field is calculated utilizing the SPH Laplacian operator \cite{Brookshaw1985} in TLSPH formalism as
\begin{equation}
	\label{eq:phasefieldLapLacian}
	\nabla^2 s_{\textbf{i}} = 2
	\sum_{\textbf{j}=1}^{N} (s_{\textbf{i}}-s_{\textbf{j}})~V_{0\textbf{j}}~ \frac{r^{s}_{0\textbf{ij}}}{|\textbf{r}_{0\textbf{ij}}|^2}~\frac{\partial W_{0\textbf{ij}}}{\partial X_{\textbf{j}}^s},
\end{equation}
where $V_{0\textbf{j}}$ is the initial volume of $\textbf{j}$, $|\textbf{r}_{0\textbf{ij}}|$ is the magnitude of the relative initial position vector $\textbf{r}_{0\textbf{ij}}=\textbf{X}_{\textbf{i}}-\textbf{X}_{\textbf{j}}$, and $r^{s}_{0\textbf{ij}}$ is the $s$-component of $\textbf{r}_{0\textbf{ij}}$. Note that Einstein summation is used for the repeated index $s$. After solving for the phase field inertia, the phase field and its first time derivative are evolved explicitly, and the results are used in updating the first Piola--Kirchhoff stress through the corresponding constitutive model. Finally, the momentum balance equation is solved for the other TLSPH field variables, e. g. $\textbf{v}$, $\dot{ \textbf{v}}$, $\textbf{u}$.

\begin{remark}
It can be observed that for soft particles ($s_{\textbf{i}}<s_l$) the deformation is assumed to be zero (\prettyref{eq:defGrad_phase}). Although in reality soft particles have obviously non-zero deformations, this is a convenient assumption that improves the numerical stability of the proposed algorithm, and does not affect the physics of the problem since the deformation gradient is only used in the stress computation which is close to zero for soft particles. From a numerical implementation perspective, soft particles stick to their neighboring undamaged particles and move with them as rigid bodies, with their only effect to the problem being through means of inertia. Obviously, a higher value of the soft particle limit would interfere with the physics of the fracture, and hence it should be selected as small as possible. Based on our experience, a value of 0.1 is enough to improve the numerical stability while preserving the physics of the problem.
\end{remark}

The overall computational procedure of the proposed framework consists of two main modules; a pre-processing module and a time integration module. In the pre-processing module, the domain is first discretized into particles, the field variables are initialized, and any preexisting discontinuities are introduced. Then, a neighbor search is performed to determine the support domain of each particle. At this stage the conventional kernel and its derivatives are computed using \prettyref{eq:TLSPHkernel}, and then the gradient correction tensor for each particle is computed from \prettyref{eq:Correction_tens}. The corrected kernel gradients are subsequently computed and substituted utilizing the $\Tilde{\nabla}W_{\textbf{ij}}=\textbf{C}_{\textbf{i}}\cdot \nabla W_{\textbf{ij}}$ relation (see \prettyref{sec:ConvSPHapproach}). The details of the prepossessing module are given in \prettyref{alg:preprocess}.
\begin{algorithm}
	\caption{Preprocessing module}
	\label{alg:preprocess}
	\begin{algorithmic}
    	\qcomnt {}
    	\subr{Preprocess}
    	\qcomnt {}
    	\pnl Get input from file or user interface
    	\pnl Discretize the continuous domain into discrete particles 
    	\pnl Introduce preexisting discontinuities 
    	\qcomnt {}
    	\For{i}{1}{Total number of particles}
    		  \qpnl Search for \textbf{j} neighbor of particle \textbf{i} in its support domain. Add \textbf{j} to the neighbor list of \textbf{i}.
    		  \qpnl Compute kernel $W_{0\textbf{ij}}$, and its gradient $\nabla_0 W_{0\textbf{ij}}$, between particles \textbf{i} and \textbf{j} using \prettyref{eq:TLSPHkernel}.
    	\EndFor
    	
    	\qcomnt {}
    	
    	\For{i}{1}{Total number of particles}
            \qpnl Compute gradient correction tensor of particle \textbf{i}, $\textbf{C}_{\textbf{i}}$, using \prettyref{eq:Correction_tens} or in index form as
            \qqcomnt {$ 
            C^{ks} _{\textbf{i}} =-\left(\sum_{\textbf{j}=1}^{N_{\textbf{i}}} r^{s}_{0\textbf{ij}} \frac{\partial W_{0\textbf{ij}}} {\partial X^{k}_{\textbf{j}}} V_{0\textbf{j}} \right)^{-1}\text{,}$} where $N_{\textbf{i}}$ is the number of neighbors of particle \textbf{i}.

    	\EndFor
    	
    	\qcomnt {}
    	\For{i}{1}{Total number of particles}
    	    \qFor{j}{1}{$N_i$}
    	        \qqpnl Correct the kernel gradient as
            $\frac{\partial W_{0\textbf{ij}}} {\partial X^{k}_{\textbf{j}}} = C^{ks} _{\textbf{i}}\frac{\partial W_{0\textbf{ij}}} {\partial X^{s}_{\textbf{j}}} $ (summation is done over index $s$)
    	    \EndqFor
    	\EndFor
    	\qcomnt {}
    	\pnl Initialize field variables
    	\pnl Apply static loads (if exist)
    	\qcomnt {}

	\end{algorithmic}
\end{algorithm}

\begin{remark}
Traditionally, in phase field for fracture, a preexisting crack is modeled either through a geometrical notch of finite width, or a prescribed value in the history functional, $\mathcal{H}_0$. Although both approaches have been extensively and successfully applied in the phase field literature, problems can potentially arise, especially in the case where the discontinuity is introduced as a physical notch. For example, in the case of meshfree and particle methods, care must be taken so that the physical discontinuity is large enough and the supports of the particles in either side of the crack do not overlap. Here, we propose a third way of modeling preexisting cracks, in which the neighbor search for particles that lie in either side of the discontinuity is restricted in that side of the crack only, as can be seen in \prettyref{fig:cracktypes}d. Simply put, we are blocking any communication among particles that lie in different sides of the preexisting crack surface $\Gamma_0$. This way the preexisting discontinuity is modeled easily and exactly without the need for a physical geometrical gap. One major advantage of this approach is that for arbitrarily shaped complex discontinuities a neighbor search restriction can be easily performed, whereas introducing geometrical gaps is fairly complicated. In this paper, we used either the physical gap or the neighbor search restriction approach and observed that both lead to similar results.

\begin{figure}[!htbp]
	\begin{center}
		\includegraphics[width=\linewidth]{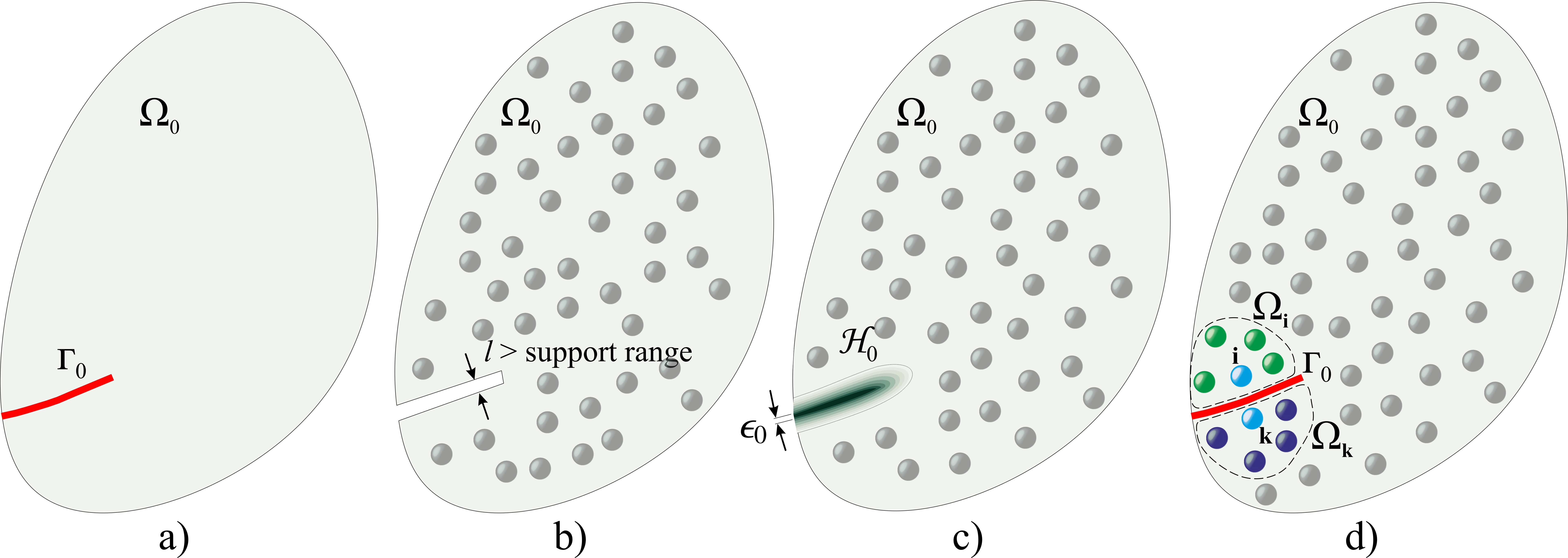}
	\end{center}
	\caption{\centering A description of a) preexisting discontinuity and its approximation by means of b) introducing a geometrical notch, c) assigning a prescribed value in the history functional, and d) restricting the neighbor search only among particles that lie in the same side of the discontinuity with the particle under consideration.}
	\label{fig:cracktypes}
\end{figure}
\end{remark}

\begin{figure}[!htbp]
	\begin{center}
		\includegraphics[width=0.7\linewidth]{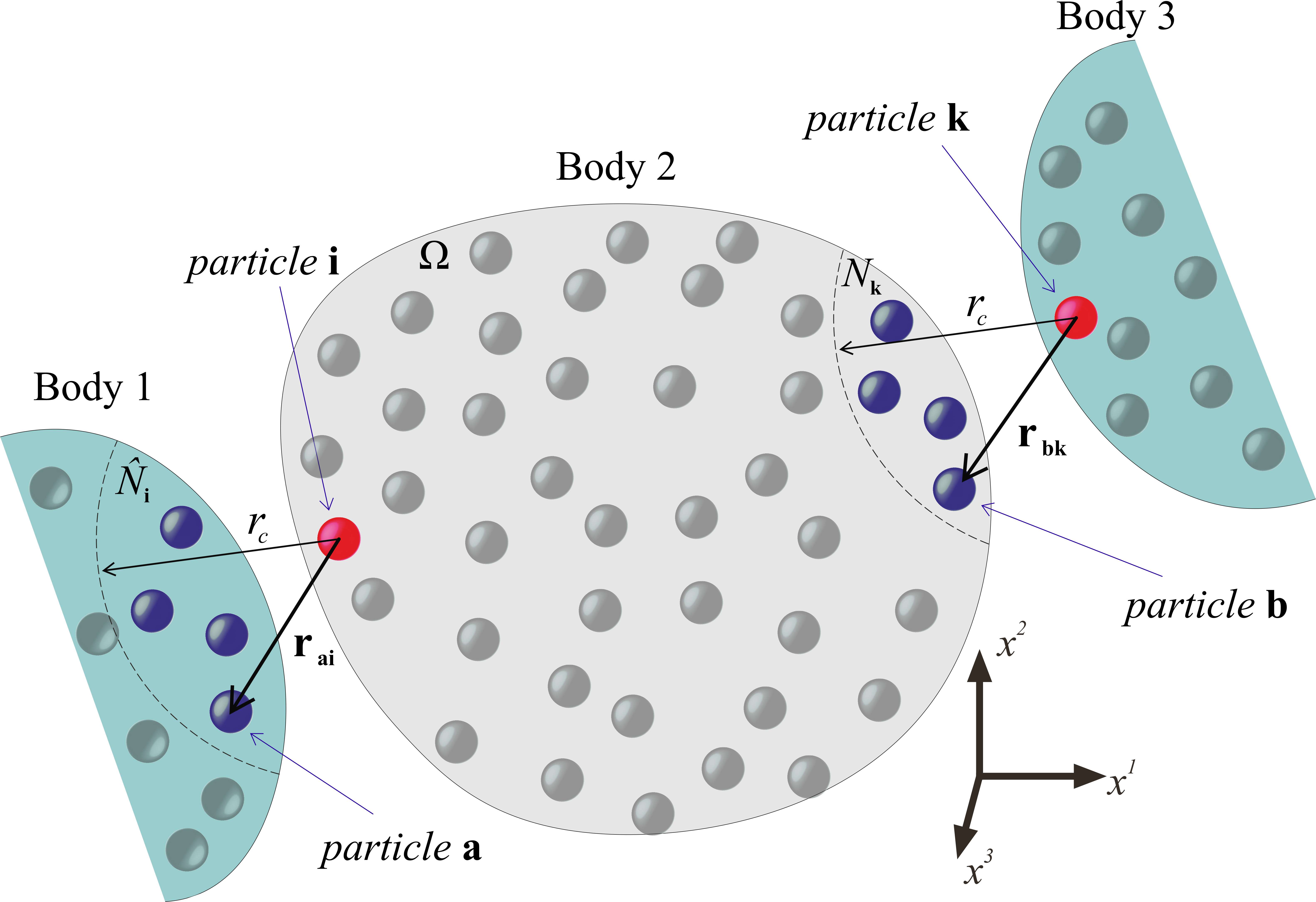}
	\end{center}
	\caption{\centering Interface details of the contact domain with multiple body contacts in SPH.}
	\label{fig:contact}
\end{figure}

When it comes to time integration, a second-order Euler predictor-corrector time integration scheme is adopted with a CFL condition of

\begin{equation}
	\Delta t\leq \text{min}\left(
	0.25\frac{h}{c_0+|\textbf{v}_{max}|},\,
	\frac{\Delta x}{c_0}
	\right),
\end{equation}
where the second term ensures the stability of the phase field solution, and $\textbf{v}_{max}$ is the maximum possible velocity of the computational domain in space and time. At the beginning of each time step the field variables are updated utilizing the predictor scheme as

\begin{equation} \label{eq:predictorv}
    \mathbf{v}_\mathbf{i}^{(t+1/2)}=\mathbf{v}_\mathbf{i}^{(t)}+0.5\Delta t\left(\frac{d\mathbf{v}_{\textbf{i}}}{dt}\right)^{(t-1/2)} \text{,}
\end{equation}
\begin{equation} 
	 \mathbf{u}_\mathbf{i}^{(t+1/2)}=\mathbf{u}_\mathbf{i}^{(t)}+0.5\Delta t\mathbf{v}_\mathbf{i}^{(t+1/2)} \text{,}
\end{equation}
\begin{equation} 
    \dot{s}_\mathbf{i}^{(t+1/2)}=\dot{s}_\mathbf{i}^{(t)}+0.5\Delta t\ddot{s}_{\textbf{i}}^{(t-1/2)} \text{,}
\end{equation}
\begin{equation} \label{eq:predictors}
	 s_\mathbf{i}^{(t+1/2)}=s_\mathbf{i}^{(t)}+0.5\Delta t\dot{s}_\mathbf{i}^{(t+1/2)} \text{,}
\end{equation}
where the superscripts $(t+1/2)$, $(t-1/2)$, and $(t)$ denote the values of the related field variable at half time step, previous half time step, and previous time step, respectively. Then, using the values of $\mathbf{u}^{(t+1/2)}$, the deformation gradient of the TLSPH particles is computed using Eqs.(\ref{eq:defGrad_phase}-\ref{eq:displacementGrad}). At this stage, an appropriate constitutive model is utilized (see \prettyref{sec:const_SVK} or \prettyref{sec:const_NH}) to compute the first Piola--Kirchhoff stress tensor as a function of the computed deformation gradient and the predicted half time step values of phase field (damage). Then, \prettyref{eq:MomentumEqTLSPHindex}, \prettyref{eq:phasefieldLapLacian}, and \prettyref{eq:phasefieldddot} are in turn solved  using the values of the variables at half-time step to compute the acceleration $(d\mathbf{v}/dt)^{(t+1/2)}$, Laplacian of the phase field $(\nabla^2s)^{(t+1/2)}$, and the phase field inertia $\ddot{s}^{(t+1/2)}$, respectively. Finally, the corrected values for the field variables are computed as
\begin{equation} \label{eq:correctorv}
    \mathbf{v}_\mathbf{i}^{(t+1)}=\mathbf{v}_\mathbf{i}^{(t)}+\Delta t\left(\frac{d\mathbf{v}_{\textbf{i}}}{dt}\right)^{(t+1/2)} \text{,}
\end{equation}
\begin{equation} 
	 \mathbf{u}_\mathbf{i}^{(t+1)}=\mathbf{u}_\mathbf{i}^{(t)}+\Delta t\mathbf{v}_\mathbf{i}^{(t+1)} \text{,}
\end{equation}
\begin{equation} 
    \dot{s}_\mathbf{i}^{(t+1)}=\dot{s}_\mathbf{i}^{(t)}+\Delta t\ddot{s}_{\textbf{i}}^{(t+1/2)} \text{,}
\end{equation}
\begin{equation} \label{eq:correctors}
	 s_\mathbf{i}^{(t+1)}=s_\mathbf{i}^{(t)}+\Delta t\dot{s}_\mathbf{i}^{(t+1)} \text{,}
\end{equation}
where the $(t+1)$ superscript denotes the current time step. Note that the initial conditions for $\dot{s}$ and $\ddot{s}$ are set to zero. \prettyref{alg:time_integration} and \prettyref{alg:const_model} present the implementation of the time integration module and the constitutive models for a single step calculation, respectively. It should be noted that for problems involving the contact of two separate bodies an additional interface force is included in the momentum equation (\prettyref{eq:MomentumEqTLSPHindex}). Hence, the momentum equation becomes
\begin{algorithm}[!htbp]
	\caption{Single step calculations of time integration module}
	\label{alg:time_integration}
	\begin{algorithmic}[1]
    	\qcomnt {}
    	\subr{Time integration}
    	\qcomnt {}
    	\pnl Apply dynamic force and/or displacement boundary conditions
    	\qcomnt {}
    	\For{i}{1}{Total number of particles}
    		  \qpnl Compute half-time step values of field variables using Eqs.(\ref{eq:predictorv}-\ref{eq:predictors})
    	\EndFor
    	\qcomnt {}
    	
    	\For{i}{1}{Total number of particles}
    	
    	   \qpnl Compute deformation gradient $\textbf{F}_{\textbf{i}}$ as in \prettyref{eq:defGrad_phase}, using half step displacement values
    	   
    	    \qpnl Utilize \prettyref{alg:const_model} to calculate $\textbf{P}_{\textbf{i}}$, $\psi_{e\textbf{i}}^{+}$, and $\psi_{e\textbf{i}}^{-}$ using the appropriate constitutive model
    	    
    	    \qIf{\psi_{e\textbf{i}}^{+}}{>}{\mathcal{H}_{\textbf{i}}}
    	        \qqpnl $\mathcal{H}_{\textbf{i}}=\psi_{e\textbf{i}}^{+}$\quad (here we update the history functional)
    	   \EndqIf
    	   
    	   \qpnl Compute $M_{\textbf{i}}$ parameter as in \prettyref{eq:Mparameter1} \quad (note that we use the upper bound value of M)
    	   
    	   \qpnl Using half-step values of the phase field parameter $s^{(t+1/2)}$ compute the Laplacian $\nabla^2 s_{\textbf{i}}$ as in \prettyref{eq:phasefieldLapLacian}
    	   
    	   \qpnl Using half step values of $s^{(t+1/2)}$ and $\dot{s}^{(t+1/2)}$ compute $\ddot{s}^{(t+1/2)}$ as in \prettyref{eq:phasefieldddot}
    	   
    	    \qpnl If contact exists, using half-step values of $\textbf{x}^{(t+1/2)}$, compute contact force for particle \textbf{i} as in \ref{app:appendixA}
    	\EndFor
    	\qcomnt {}
    	
    	\For{i}{1}{Total number of particles}
    	    \qpnl Using half step values of field variables compute acceleration $(d\mathbf{v}_{\textbf{i}}/dt)^{(t+1/2)}$ as in \prettyref{eq:MomentumEqTLSPHindex} (or \prettyref{eq:momentumcontact} for problems involving contact interactions)
    	    \qpnl Correct field variables using Eqs.(\ref{eq:correctorv}-\ref{eq:correctors})
    	\EndFor
    	\qcomnt{}
    	\pnl Write results (occasionally)
    	
	\end{algorithmic}
\end{algorithm}
\begin{algorithm}[!htbp]
	\caption{Implementation of the Saint-Venant Kirchhoff and neo-Hookean hyperelastic constitutive models}
	\label{alg:const_model}
	\begin{algorithmic}[1]
    	\qcomnt {}
    	\subr{If constitutive model is Saint-Venant Kirchhoff:}
        	\qpnl Compute Green--Lagrange strain tensor as $ 
            E_{\mathbf{i}}^{ks} =\frac{1}{2}\left(
            F_{\mathbf{i}}^{mk}F_{\mathbf{i}}^{ms}+\delta^{ks}
            \right)$ (summation is done over index m)
        	    
    	    \qpnl Decompose $\textbf{E}_{\textbf{i}}$ into $\textbf{E}_{\textbf{i}}^{+}$ and $\textbf{E}_{\textbf{i}}^{-}$ using Eqs.(\ref{eq:GLstratinP}-\ref{eq:GLstratinN})
    	    
    	    \qpnl Compute $\psi_{e\textbf{i}}^{+}$ and $\psi_{e\textbf{i}}^{-}$ using Eqs.(\ref{eq:elasEnergyP}-\ref{eq:elasEnergyN})
    	    
    	    \qpnl Compute $\textbf{S}_{\textbf{i}}^{+}$ and $\textbf{S}_{\textbf{i}}^{-}$ using \prettyref{eq:2PK_SVK1}
    	    
    	    \qpnl Compute $\textbf{S}_{\textbf{i}}= \left(s_{\textbf{i}}^{(t+1/2)}\right)^2 \textbf{S}_{\textbf{i}}^{+} + \textbf{S}_{\textbf{i}}^{-}$
    	    
    	    \qpnl Compute $P^{ks}_{\textbf{i}}=F^{km}_{\textbf{i}}S^{ms}_{\textbf{i}}$ (summation is done over index m)
	        \qcomnt {}
	    \subr{Else if constitutive model is neo-Hookean:}
	        \qpnl Compute $J$, $\textbf{b}$, and $\overline{\textbf{b}}$ using \prettyref{eq:jacobian}, \prettyref{eq:bmatrix}, and \prettyref{eq:bbarmatrix}, respectively.
	        
	        \qpnl Compute $U(J)$ and $\overline{\psi}_e(\overline{\textbf{b}})$ using \prettyref{eq:UjNeoHok} and \prettyref{eq:ksiNeoHok}, respectively.
	        
	         \qpnl Compute $\psi_{e\textbf{i}}^{+}$ and $\psi_{e\textbf{i}}^{-}$ using \prettyref{eq:ksipNH} and \prettyref{eq:ksinNH}, respectively. 
	         
            \qpnl Compute $\textbf{S}_{\textbf{i}}$ from \prettyref{eq:2PK_NH2} using half-step values of phase field, $s_{\textbf{i}}^{(t+1/2)}$.
            
            \qpnl Compute $P^{ks}_{\textbf{i}}=F^{km}_{\textbf{i}}S^{ms}_{\textbf{i}}$ (summation is done over index m)
	    \qcomnt {}
	    \subr{Else:}
    	    \qpnl Apply a different constitutive model here to calculate $\textbf{P}_{\textbf{i}}$, $\psi_{e\textbf{i}}^{+}$, and $\psi_{e\textbf{i}}^{-}$

	\end{algorithmic}
\end{algorithm}

\begin{equation} \label{eq:momentumcontact}
	\frac{d v^k_\textbf{i}}{dt} = \sum_{\textbf{j}=1}^{N_{\textbf{i}}} 
	m_{0\textbf{j}}
	\left( 
	\frac{P^{ks}_\textbf{i}}{\rho^2_{0\textbf{i}}} +
	\frac{P^{ks}_\textbf{j}}{\rho^2_{0\textbf{j}}} +
	P^{ks}_{v\textbf{ij}}
	\right)
	~\frac{\partial W_{0\textbf{ij}}}{\partial X_{\textbf{j}}^s} + b_{0\textbf{i}}^k + \sum_{\textbf{a}=1}^{\hat{N}_{\textbf{i}}} f_{\textbf{ai}} r_{\textbf{ai}}^k~
\end{equation}
where $\hat{N}_{\textbf{i}}$ is the total number of particles located at a separate body within the contact distance, $r_c$, of particle $\textbf{i}$, $r_{\textbf{ai}}^k$ is the k component of the relative position vector $\textbf{r}_{\textbf{ai}}=\textbf{x}_{\textbf{a}}-\textbf{x}_{\textbf{i}}$ between particles $\textbf{a}$ and $\textbf{i}$ calculated using the Eulerian coordinates $\textbf{x}$, as given in  \prettyref{fig:contact}, and $f_{\textbf{ai}}$ is a scalar multiplier whose details are given in \ref{app:appendixA}.  For the cases involving contact interaction in this paper we set $r_c = 2.5 \Delta x$ unless stated otherwise.

\section{Numerical examples}
\label{sec:results}
In this section, we first apply the proposed framework to model various challenging problems from the literature. We then illustrate further capabilities of our model by addressing the failure occurrences in a three dimensional impact problem. All problems except the last one are simulated in two dimensions under plane strain assumption with a CPU parallel code. The simulations were carried out on a 40-core node on SeaWulf cluster located at the Institute for Advanced Computational Science (IACS) at Stony Brook University. The animated videos of the simulations are provided in Electronic
Annex I (\ref{app:appendixB}).

\subsection{Symmetric three point bending test}

\begin{figure}[!htbp]
	\begin{center}
		\includegraphics[width=0.5\linewidth]{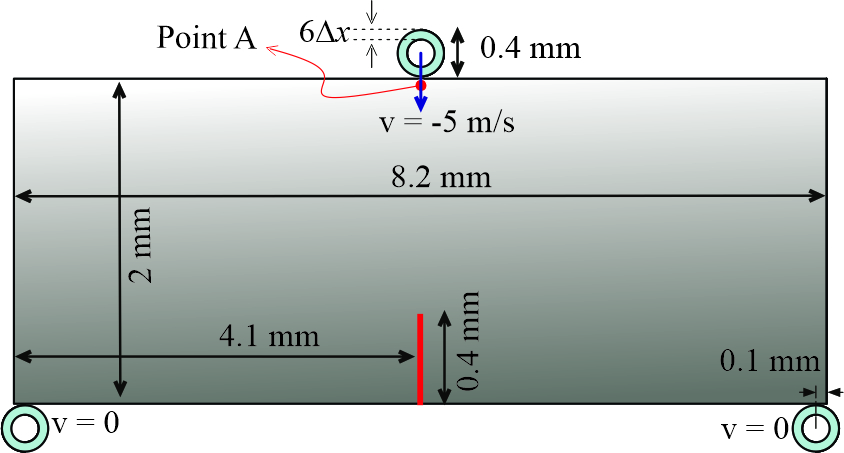}
	\end{center}
	\caption{\centering Symmetric three point bending test. Geometry and loading conditions.}
	\label{fig:3PbendingGeosym}
\end{figure}

\begin{figure}[!htbp]
	\begin{center}
		\includegraphics[width=\linewidth]{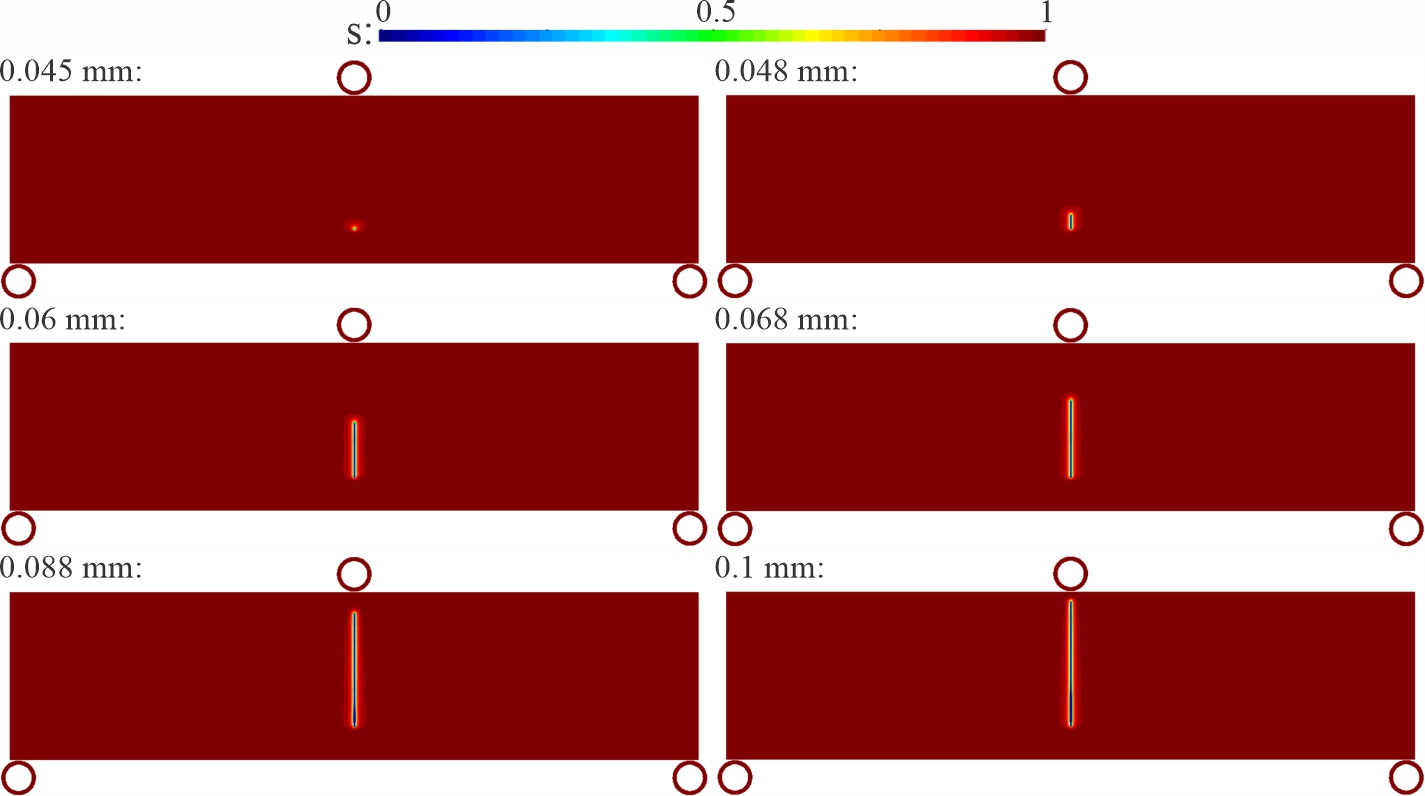}
	\end{center}
	\caption{\centering Symmetric three point bending test. Contours of the phase field parameter (crack path) for different values of the $X^2$-displacement of Point A. The snapshots correspond to computations performed with the Saint-Venant Kirchhoff constitutive model.}
	\label{fig:3Pbendingsymphasefield}
\end{figure}

\begin{figure}[!htbp]
	\begin{center}
		\includegraphics[width=0.5\linewidth]{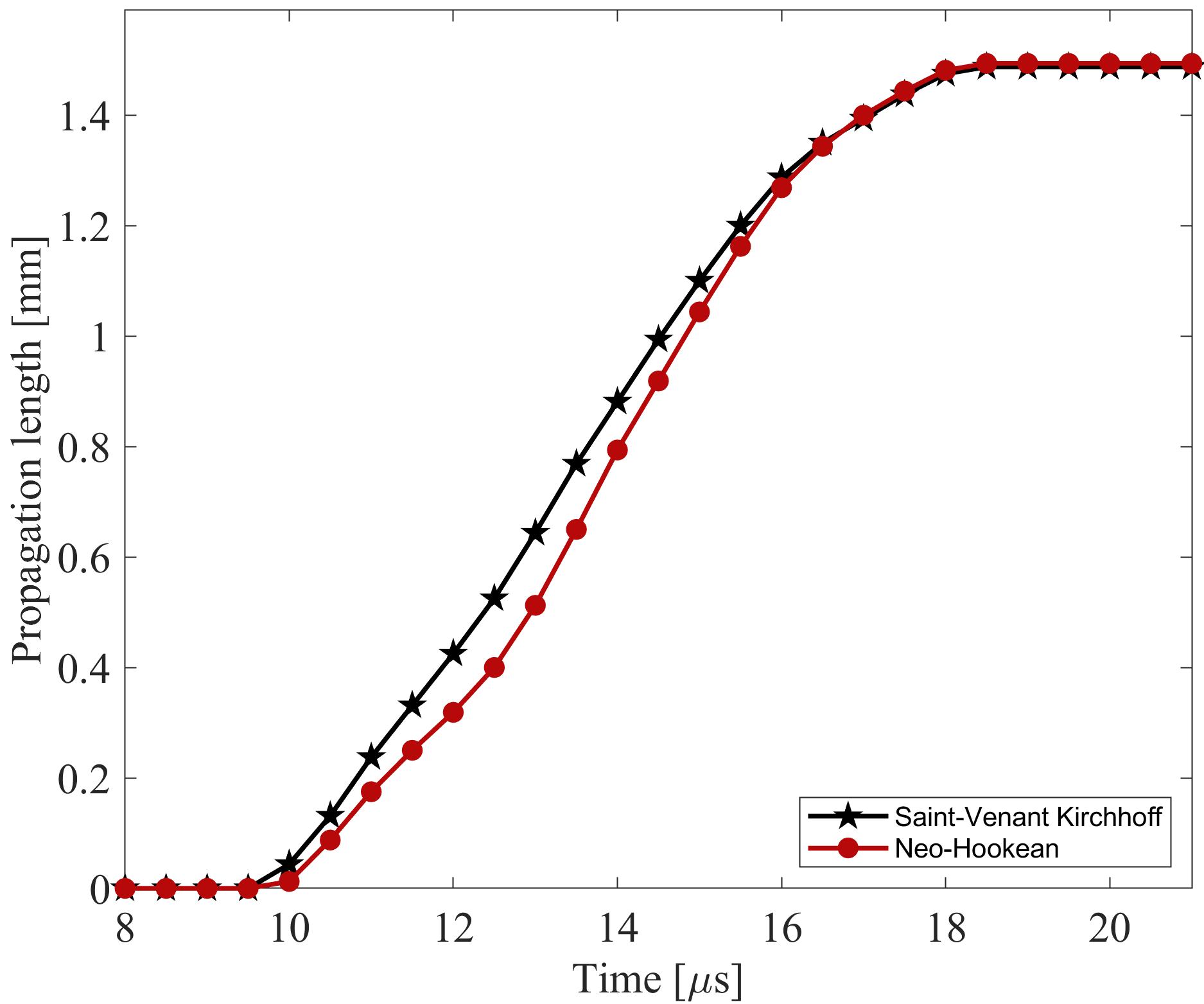}
	\end{center}
	\caption{\centering Symmetric three point bending test. Evolution of the propagation length of the crack for the Saint-Venant Kirchhoff and Neo-Hookean constitutive models.}
	\label{fig:3Pbendingsymgraph}
\end{figure}

We first solve the crack propagation in a plate subjected to symmetric three point bending \cite{ambati2015review,miehe2010phase}. The length and width of the plate are set to $L_x=8.2$~mm and $L_y=2$~mm, respectively. The material parameters $\lambda = 12$~GPa, $\mu=8$~GPa, density $\rho_0 = 1190$~kg/m$^3$, and $G_c=0.54$~kJ/m$^2$ are adopted from \cite{miehe2010phase}. The contact parameters, length scale, and time step values, are set to $K_p=9\e{12}$ and $r_0=2.5\Delta x$, $\epsilon_0=\Delta x$, and $\Delta t=0.1$~ns, respectively. A preexisting crack with a length of $0.4$~mm is located at the lower central portion of the plate and is modeled by restricting the neighbor search of the particles around the discontinuity region, as suggested in \prettyref{sec:procedure} and \prettyref{fig:cracktypes}d. We enforce the boundary conditions through contact with rigid circular bands with outer diameters of $0.4$~mm at three regions, as shown in \prettyref{fig:3PbendingGeosym}. The contact force is applied in the $X^2$-direction only, and a velocity of $-5$~m/s in the $X^2$-direction is assigned to the upper band, whereas the two lower bands are fixed in the $X^1$- and $X^2$-directions.
We set the artificial viscosity parameters (see \prettyref{sec:GovEqTLSPH}) to $\beta_1=0.2$ and $\beta_2=0$, and a phase field limit of $s_l=0.1$ is chosen for soft particles. The domain is discretized into $423,665$ particles with a particle spacing of $\Delta x= 6.25$~$\mu$m, out of which $3,825$ particles belong to the contact bodies. Both the Saint-Venant Kirchhoff and Neo-Hookean constitutive models are employed. \prettyref{fig:3Pbendingsymphasefield} shows the snapshots of the phase field parameter at different stages of the propagation.
As can be seen, the propagation starts from the preexisting crack tip, continues in a vertical direction, and stops near the upper boundary of the plate. The start and end time of the propagation correspond to a displacement of $0.045$~mm and $0.1$~mm of Point A (see \prettyref{fig:3PbendingGeosym}) in the $X^2$-direction, which agrees well with the ones reported in \cite{miehe2010phase}. In \prettyref{fig:3Pbendingsymgraph} we compare the propagation length predicted by Saint-Venant Kirchhoff and Neo-Hookean constitutive models, in which both models deliver similar propagation curves.

\subsection{Asymmetric three point bending test}

\begin{figure}[!htbp]
	\begin{center}
		\includegraphics[width=0.5\linewidth]{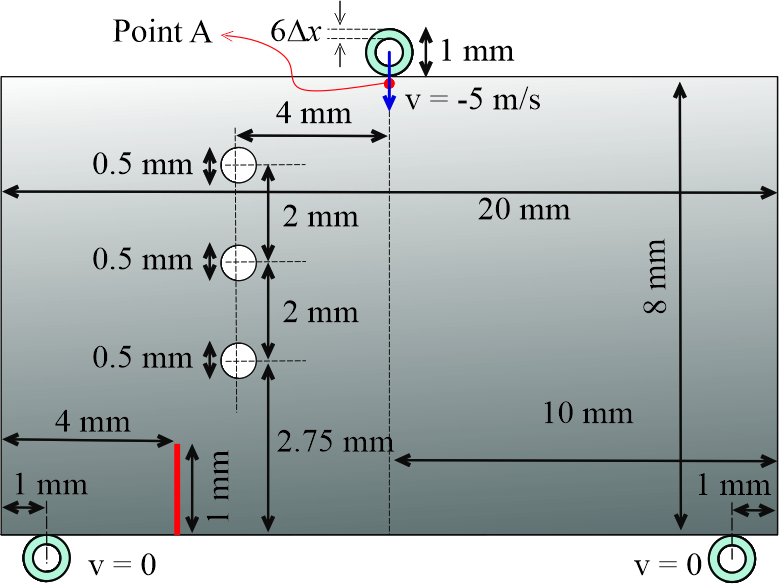}
	\end{center}
	\caption{\centering Asymmetric three point bending test. Geometry and loading conditions.}
	\label{fig:asym3PbendingGeo}
\end{figure}

\begin{figure}[!htbp]
	\begin{center}
		\includegraphics[width=0.85\linewidth]{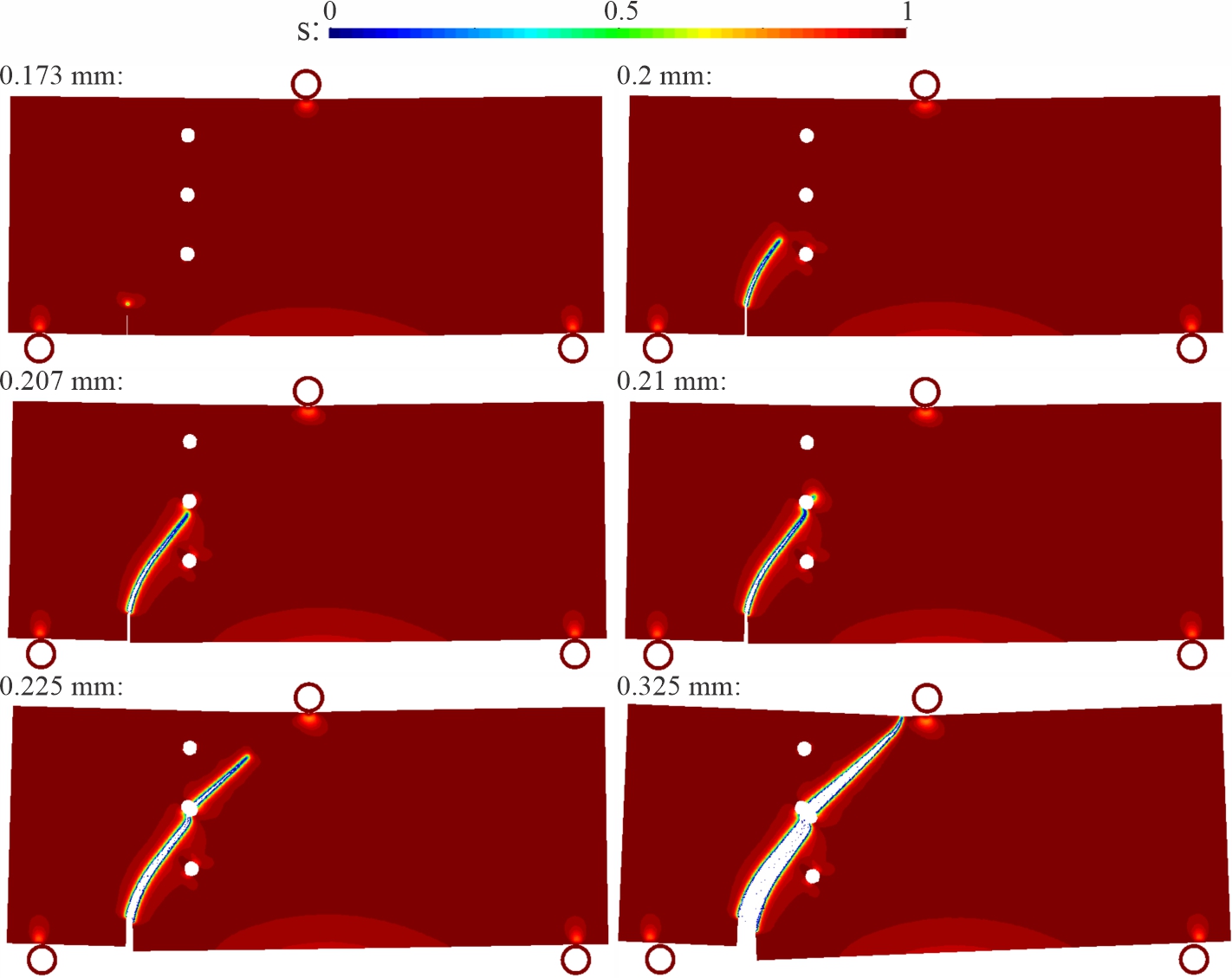}
	\end{center}
	\caption{\centering Asymmetric three point bending test. Contours of the phase field parameter (crack path) for different values of the $X^2$-displacement of Point A.}
	\label{fig:3Pbendingasymphasefield}
\end{figure}

\begin{figure}[!htbp]
	\begin{center}
		\includegraphics[width=0.4\linewidth]{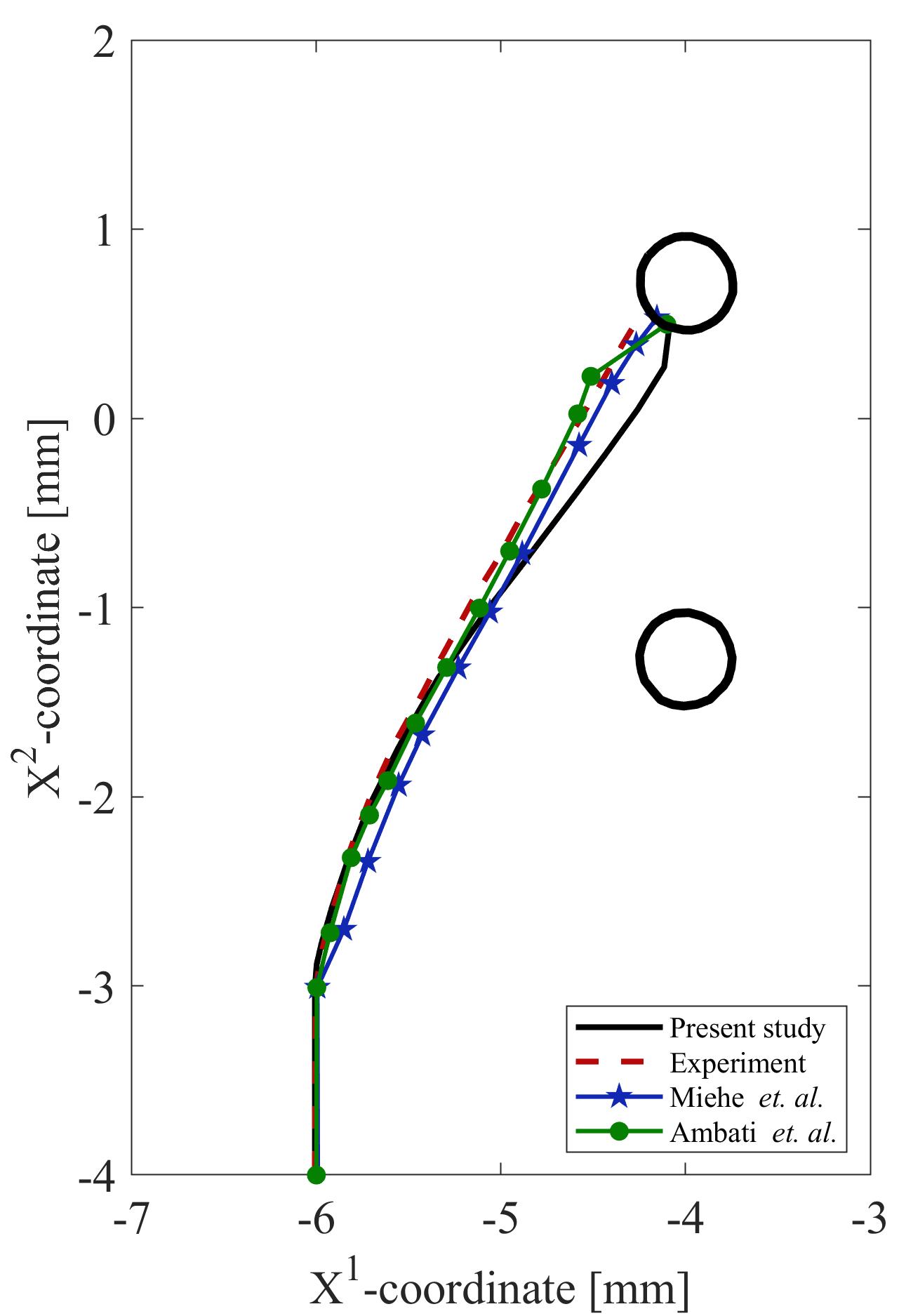}
	\end{center}
	\caption{\centering Asymmetric three point bending test. Crack path compared with experimental \cite{BITTENCOURT1996} and phase field results of Miehe \etal\cite{miehe2010phase} and Ambati \etal\cite{ambati2015review}.}
	\label{fig:3Pbendingasymgraph}
\end{figure}

Here, we model the asymmetric three point bending test \cite{BITTENCOURT1996,ambati2015review,miehe2010phase}. The boundary conditions are similar to the previous case and are given in \prettyref{fig:asym3PbendingGeo}. The geometric properties and material parameters are adopted as $L_x=20$~mm, $L_y=8$~mm, $\lambda = 12$~GPa, $\mu=8$~GPa, $\rho_0 = 1190$~kg/m$^3$, and $G_c=1$~kJ/m$^2$. 
We model the preexisting crack by restricting the neighbor search of the particles around the discontinuity region. The contact force is applied in the $X^2$-direction with contact parameters of $K_p=12\e{12}$ and $r_0=2.5\Delta x$. Other simulation parameters are set to $\epsilon_0=\Delta x$, $\Delta t=0.1~ns$, $\beta_1=0.2$, $\beta_2=0$, and $s_l=0.1$. We discretize the domain into approximately 400K particles with a particle spacing of $\Delta x = 0.02$ mm. The Saint-Venant Kirchhoff constitutive model is employed. \prettyref{fig:3Pbendingasymphasefield} shows the snapshots of the phase field for the deformed plate at different stages of the propagation.
As can be seen, the crack initiation occurs when point A has displaced by approximately $0.17$~mm in the $X^2$-direction.  Then, the crack joins the middle hole after a displacement of approximately $0.21$~mm, and the complete rupture is recorded at a displacement of around $0.325$~mm. The crack first propagates on a somewhat curved path from the preexisting crack tip to the middle hole, and then on an almost linear path until the full rupture. \prettyref{fig:3Pbendingasymgraph} compares the experimentally \cite{BITTENCOURT1996} and numerically \cite{ambati2015review,miehe2010phase} obtained crack paths with the one predicted by the present approach. 
As can be seen, the predicted crack path agrees very well with both the experimental and numerical studies found in the literature. There is a small deviation observed within a small region near the middle hole. One potential reason for this could be the imperfect representation of the circular holes in the discretized particle domain, which can alter the local strain energy concentrations. Nevertheless, the recorded deviation is negligible and does not change the overall physical behavior of the crack.

\subsection{Dynamic crack branching}
\begin{figure}[!htbp]
	\begin{center}
		\includegraphics[width=0.5\linewidth]{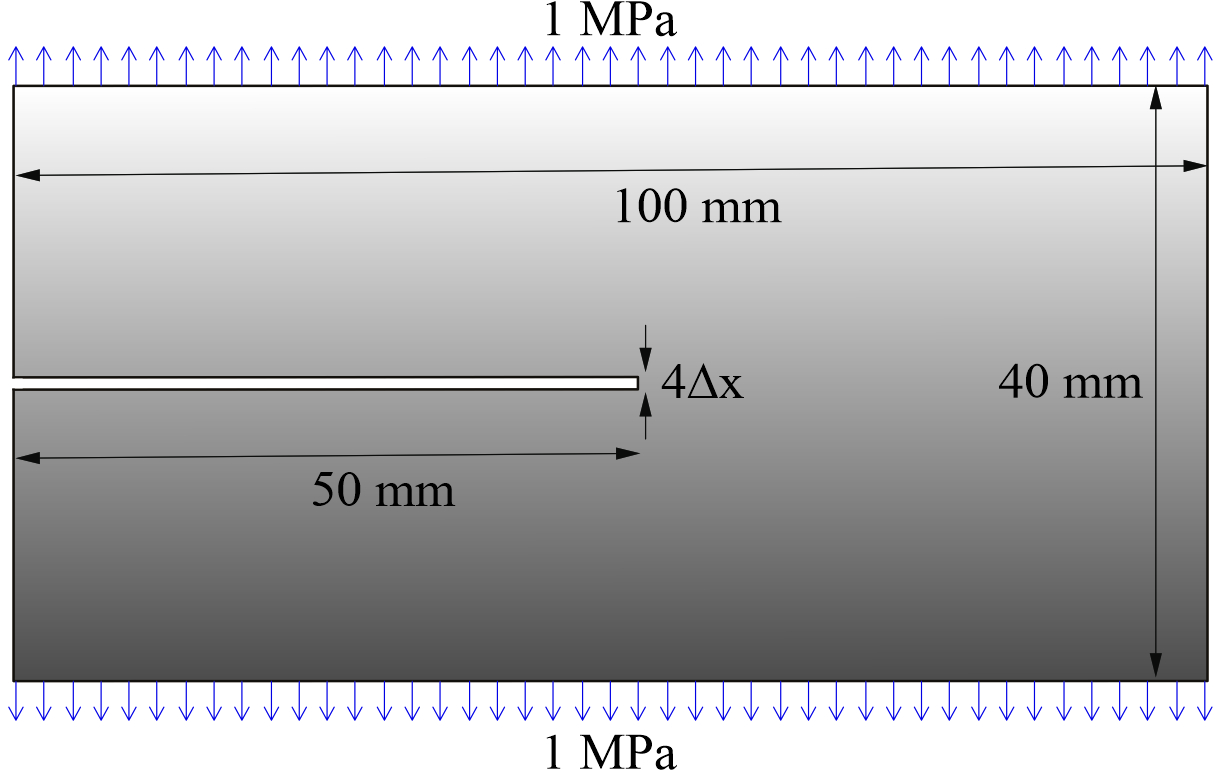}
	\end{center}
	\caption{\centering Dynamic crack branching. Geometry and loading conditions.}
	\label{fig:crackBranchGeo}
\end{figure}

The dynamic crack branching problem has been widely investigated for different boundary conditions and loading setups \cite{Borden2012,Kamensky2018,moutsanidis2018hyperbolic,Rahimi2020,Rahimi2020b}. In the present effort we focus on the setup used in \cite{Borden2012,Kamensky2018,moutsanidis2018hyperbolic}, where a plate with a notch is subjected to a tensile surface loading. We set the material properties of the plate to $E=32$~GPa, $\rho = 2450$~kg/m$^3$, $\nu = 0.2$, and $G_c=3$~J/m$^2$, whereas the Saint-Venant Kirchhoff constitutive model is employed. The length and width of the plate are taken as $L_x=100$~mm and $L_y=40$~mm, respectively. The preexisting notch is introduced as a geometrical discontinuity with a width of $4\Delta x$ and length of $50$~mm located at the center. The plate is assumed to be under a tensile load of $1$~MPa applied to its upper and lower surfaces, as depicted in \prettyref{fig:crackBranchGeo}. The length scale and time step values are set to $\epsilon_0=0.09375$~mm and $\Delta t=2.5$~ns, respectively. We set the artificial viscosity parameters (see \prettyref{sec:GovEqTLSPH}) to $\beta_1=0.04$ and $\beta_2=0$, and a phase field limit of $s_l=0.1$ is chosen for soft particles.
To demonstrate the resolution independence of the crack pattern in the phase field method, we consider four different particle resolutions\footnote{Care should be taken when choosing the minimal resolution and time step. Although the crack path in phase field is independent of the resolution, extremely low particle resolutions may alter the path of the crack because the fracture energy is regularized over a larger length scale (cracks are wider), which can in turn affect the physics of the problem. Additionally, in TLSPH, for lower resolutions and large time step values, higher amount of artificial viscosity is needed to overcome the zero-energy mode discrepancy \cite{monaghan1983374}, thus, leading to a higher amount of nonphysical forces in the domain. This can then alter the original speed of the crack propagation, and in some cases even the path.}, R1-R4, with initial particle spacing of $\Delta x=0.1$~mm, $\Delta x=0.08$~mm, $\Delta x=0.0625$~mm, and $\Delta x=0.05$~mm, leading to a discretization of the domain into $398,000$, $622,500$, $1,020,800$, and $1,596,000$ particles, respectively.
\begin{figure}[!htbp]
	\begin{center}
		\includegraphics[width=0.8\linewidth]{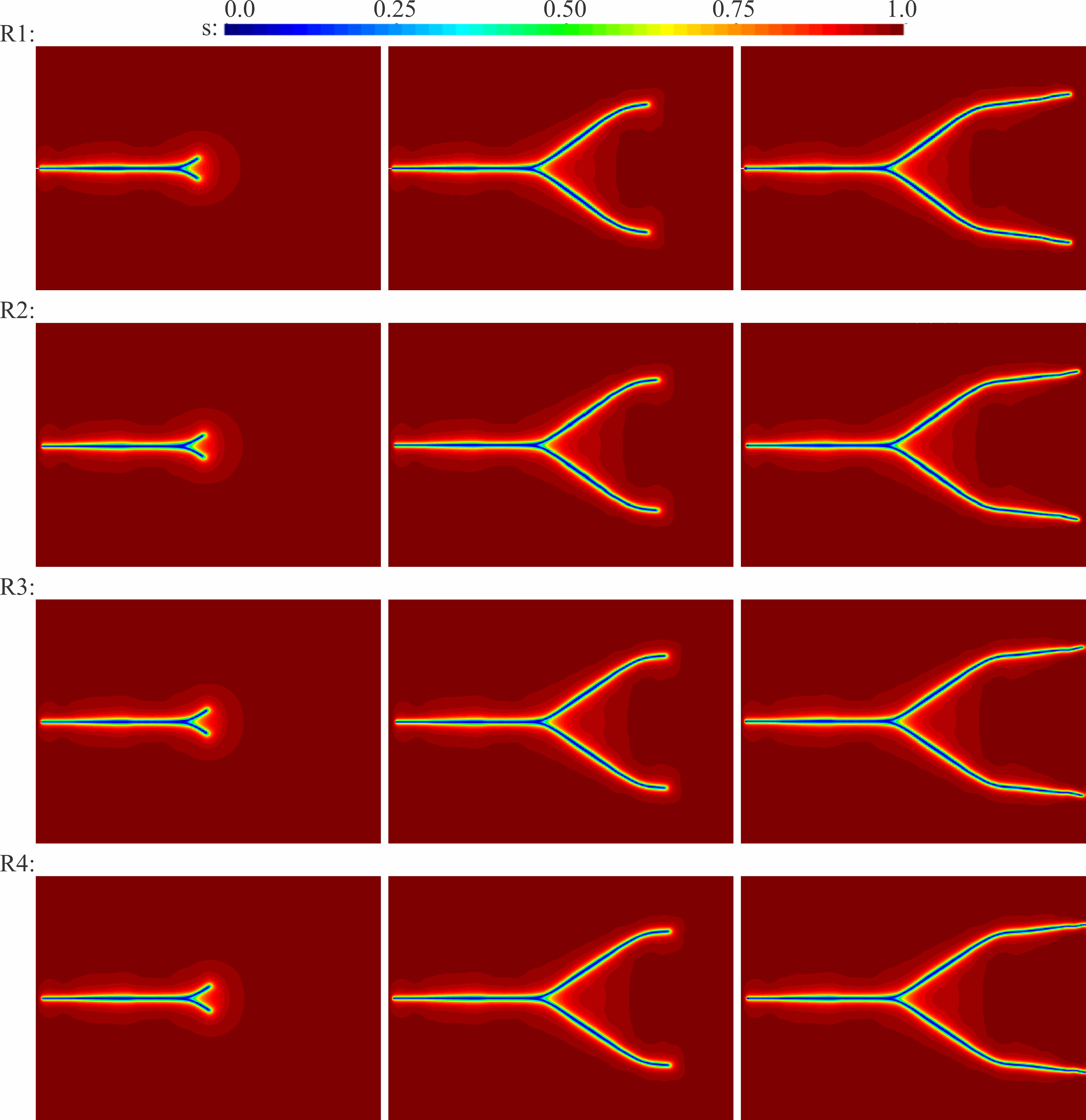}
	\end{center}
	\caption{\centering Dynamic crack branching. Contours of the phase field parameter (crack path) at $55~\mu$s (left), $80~\mu$s (middle), and $105~\mu$s (right), for R1-R4 resolutions. The snapshots are presented for the right half of the plate.}
	\label{fig:crackBranchingphasefield}
\end{figure}

\begin{figure}[!htbp]
	\begin{center}
		\includegraphics[width=0.8\linewidth]{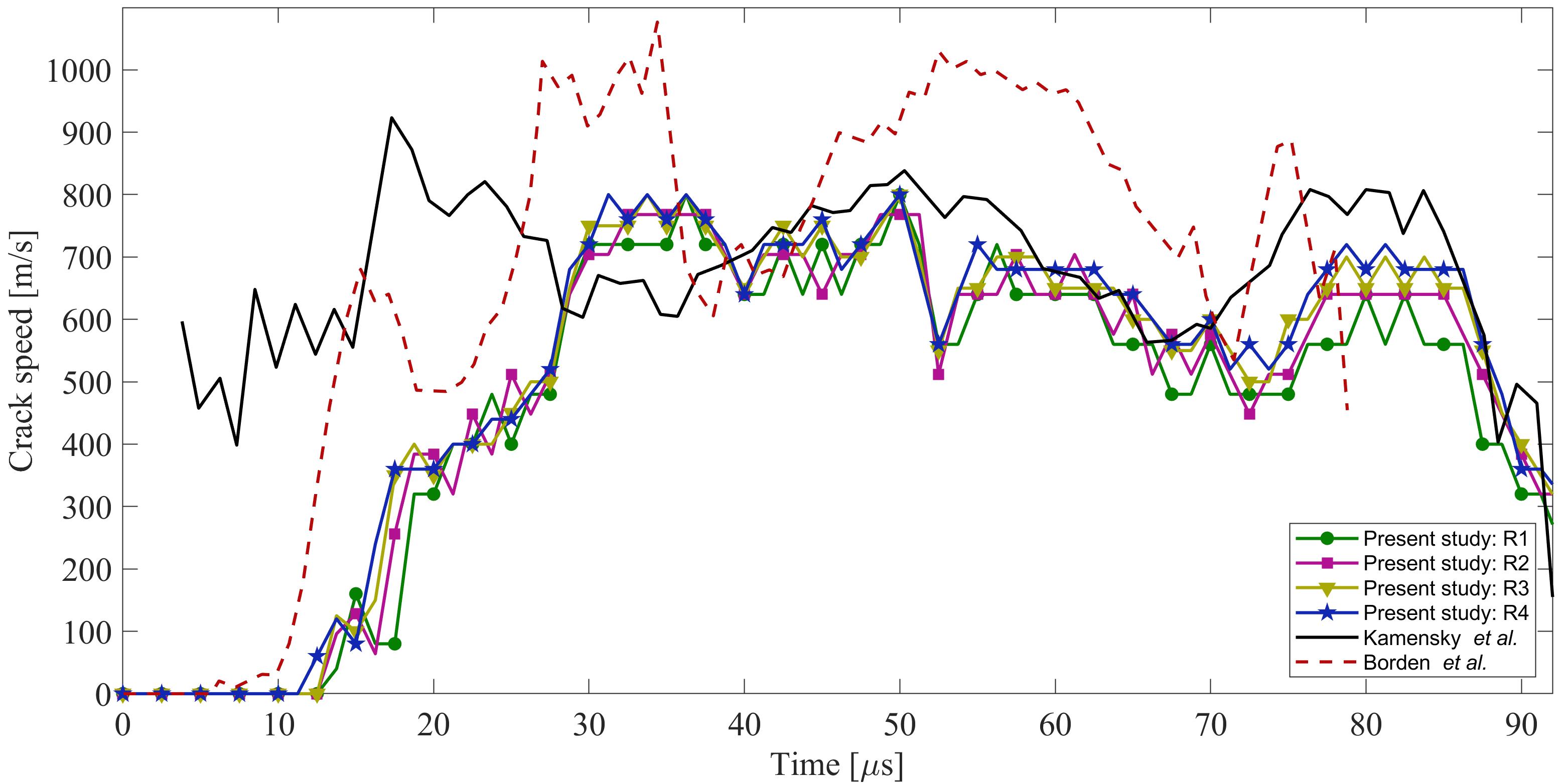}
	\end{center}
	\caption{\centering Dynamic crack branching. $X^1$- component of crack velocity for resolutions R1-R4 compared with the phase field results of Kamensky \etal\cite{Kamensky2018} and Borden \etal\cite{Borden2012}.}
	\label{fig:crackBranchGraph}
\end{figure}

\prettyref{fig:crackBranchingphasefield} presents several snapshots of the phase field contours for R1-R4 resolutions at different stages of the crack propagation for the right half of the plate (the left half stays intact and we therefore omit that region). The crack propagation starts from the tip of the notch at approximately $t=10$~$\mu$s and continues in a horizontal manner. At approximately $t=50$~$\mu$s, branching occurs, and the fracture propagates at an angle of $\approx 37^o$. At later stages, the propagation angles of the two branches tend to change back to horizontal. As expected, an identical crack pattern is predicted for all different resolutions, highlighting phase field's feature of mesh independence, which is hard to achieve when local damage models are employed. At the same time, the predicted crack patterns are in excellent agreement with the ones reported in \cite{Borden2012,Kamensky2018,moutsanidis2018hyperbolic}. However, the crack in resolutions R1 and R2 fails to reach the end of the plate in the given time. This is a common behavior seen for low discretization resolutions and it is attributed to the fact that for lower resolutions the fracture zone is slightly wider and hence there is more fracture energy per unit length.
Similar behavior was also reported in \cite{Borden2012,Kamensky2018,moutsanidis2018hyperbolic}. \prettyref{fig:crackBranchGraph} compares the crack propagation speed predicted by the current model with the ones of \cite{Borden2012,Kamensky2018}. A close agreement is observed between the present results and the ones of \cite{Kamensky2018}. In \cite{Borden2012}, however, a slightly higher propagation speed is reported. This difference is attributed to the absence of the inertia term ($\ddot{s}$) in the phase-field PDE utilized in \cite{Borden2012} (see \prettyref{sec:PhaseField}). Obviously, the inclusion of the inertia term in the current effort (i.e. the use of a hyperbolic phase field PDE instead of an elliptic one, as described in Section \ref{sec:PhaseField}) introduces an extra amount of work to be carried out in order for the phase field to propagate, thereby causing a decrease in the propagation speed. However, it should be pointed out that no experimental data exist for this problem, therefore it is not clear which model’s crack propagation speed is more accurate. Nonetheless, as mentioned in \cite{Kamensky2018}, the elliptic models tend to overestimate the crack propagation speeds because they ignore the rate toughening effect in brittle fracture.
It is also worth mentioning that the strategy employed to track the crack tip, and hence the propagation speed, highly affects the recorded crack initiation time. Here, we consider the furthest particle in the $+X^1$-direction with a phase field value of $s<0.7$ to be the tip of the crack. Therefore, as soon as the phase field value decreases to 0.7, the propagation is assumed to initiate. However, the choice of this value does not affect the computed propagation speed since the speed is computed between the two iso-curves of the same value. 

\begin{figure}[!htbp]
	\begin{center}
		\includegraphics[width=0.6\linewidth]{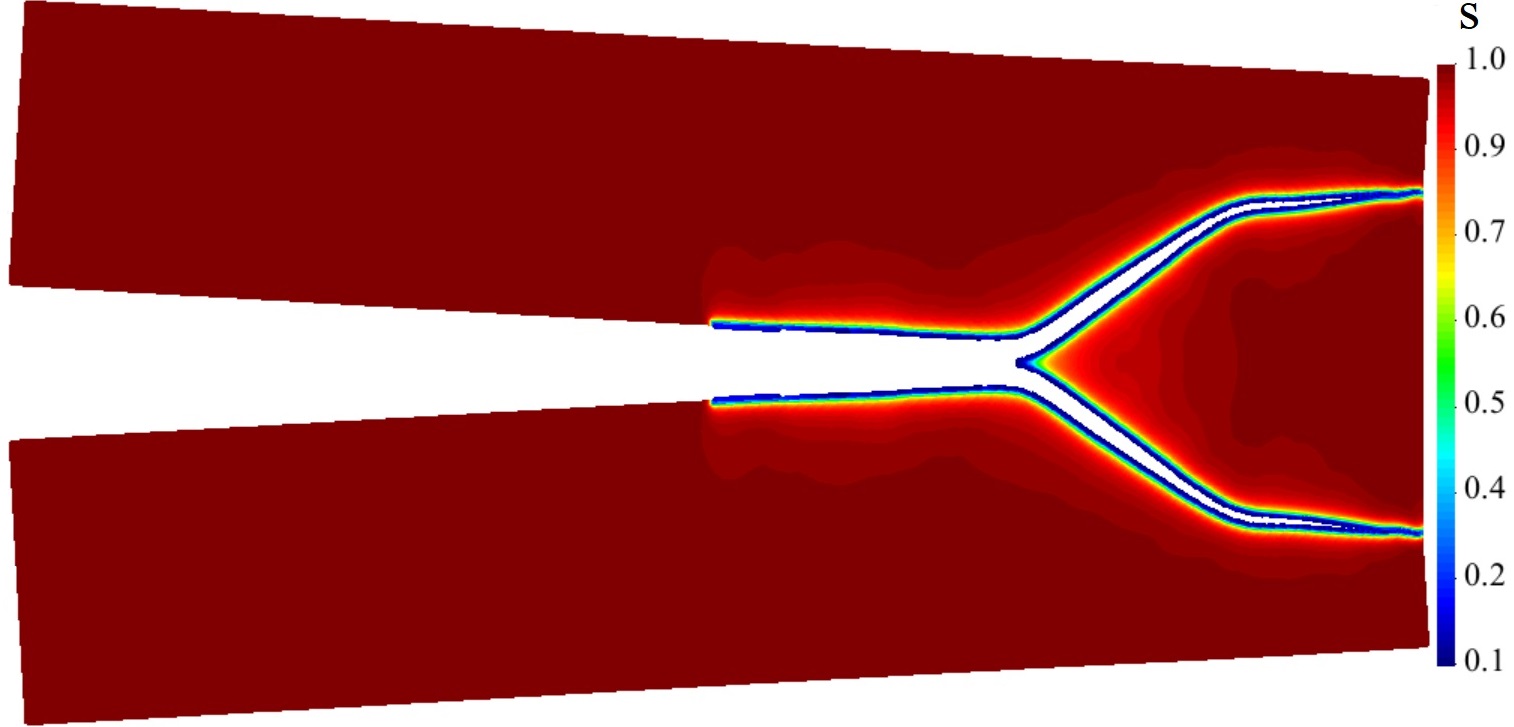}
	\end{center}
	\caption{\centering Dynamic crack branching. Opening of the crack and phase field distribution at the failure moment ($\approx 105$~$\mu$s) for resolution R4. Particles with $s<s_l=0.1$ are not shown. The deformation is magnified by a factor of $50$.}
	\label{fig:crackBranchcrackopening}
\end{figure}

In \prettyref{fig:crackBranchcrackopening} the post-processed final results of the phase field and the crack opening are given for resolution R4, where the particles with $s<s_l=0.1$ are not shown, and the deformation is scaled by a factor of 50. As can be seen, after the complete failure, the domain separates into three regions, marking the end of the propagation.

\subsection{Kalthoff--Winkler experiment}

\begin{figure}[!htbp]
	\begin{center}
		\includegraphics[width=0.45\linewidth]{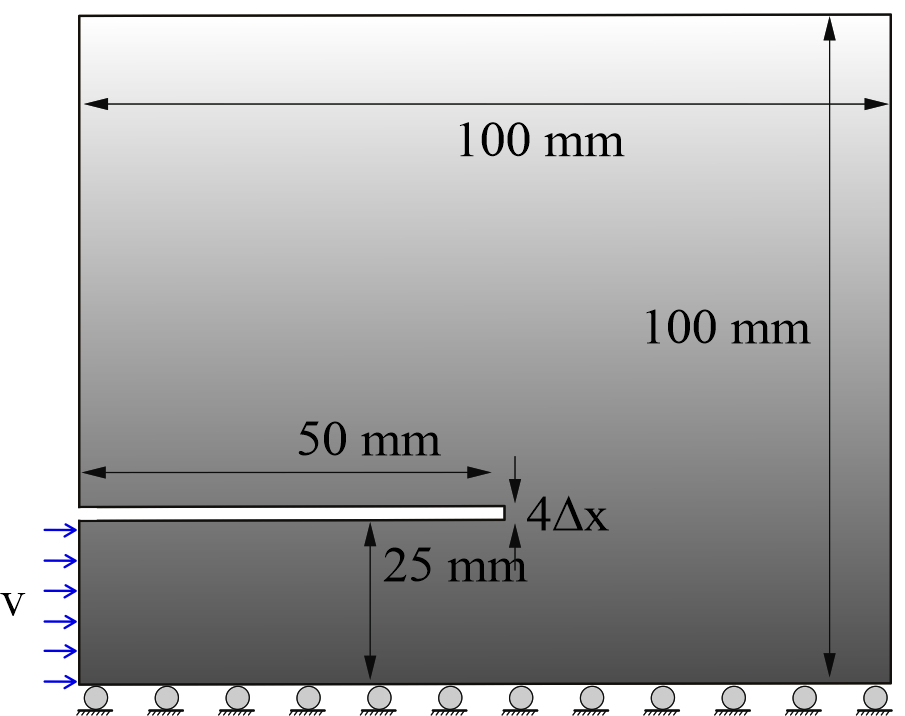}
	\end{center}
	\caption{\centering Kalthoff--Winkler experiment. Geometry and loading conditions.}
	\label{fig:kalthoffwinkGeo}
\end{figure}

\begin{figure}[!htbp]
	\begin{subfigure}{\linewidth}
		\begin{center}
			\includegraphics[width=0.7\linewidth]{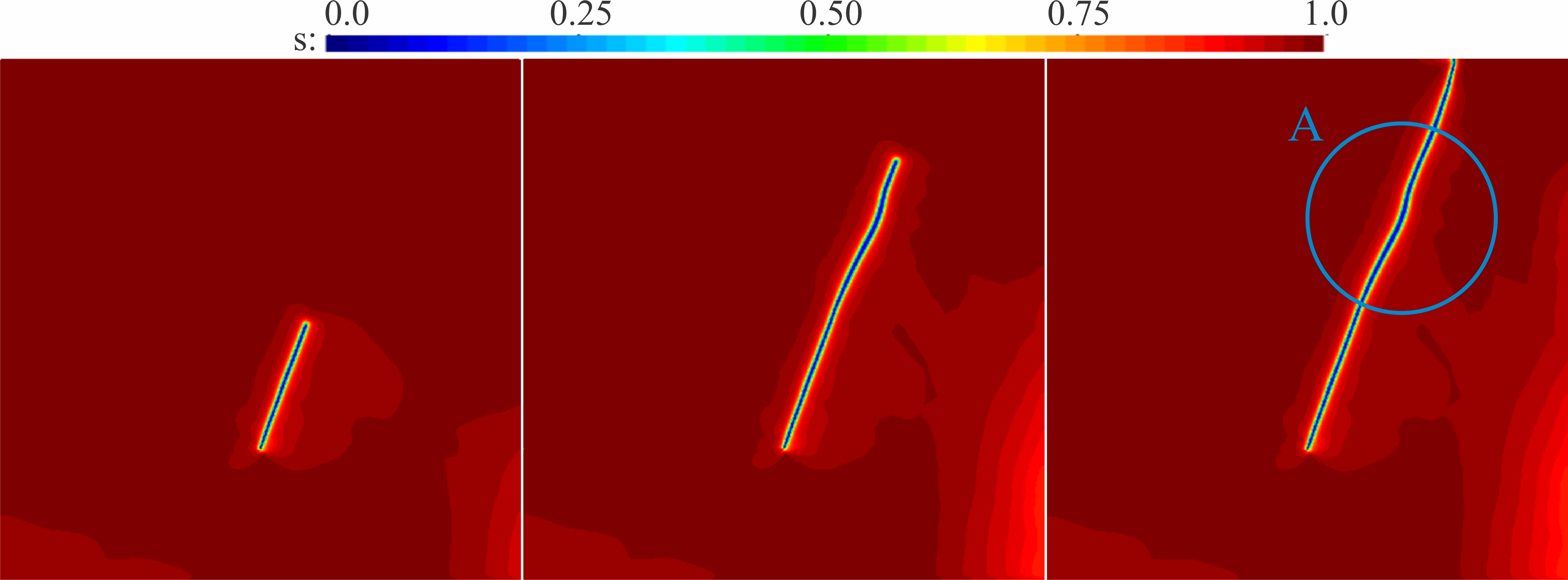}
		\end{center}
		\caption{ }
	\end{subfigure}
	\begin{subfigure}{\linewidth}
		\begin{center}
			\includegraphics[width=0.7\linewidth]{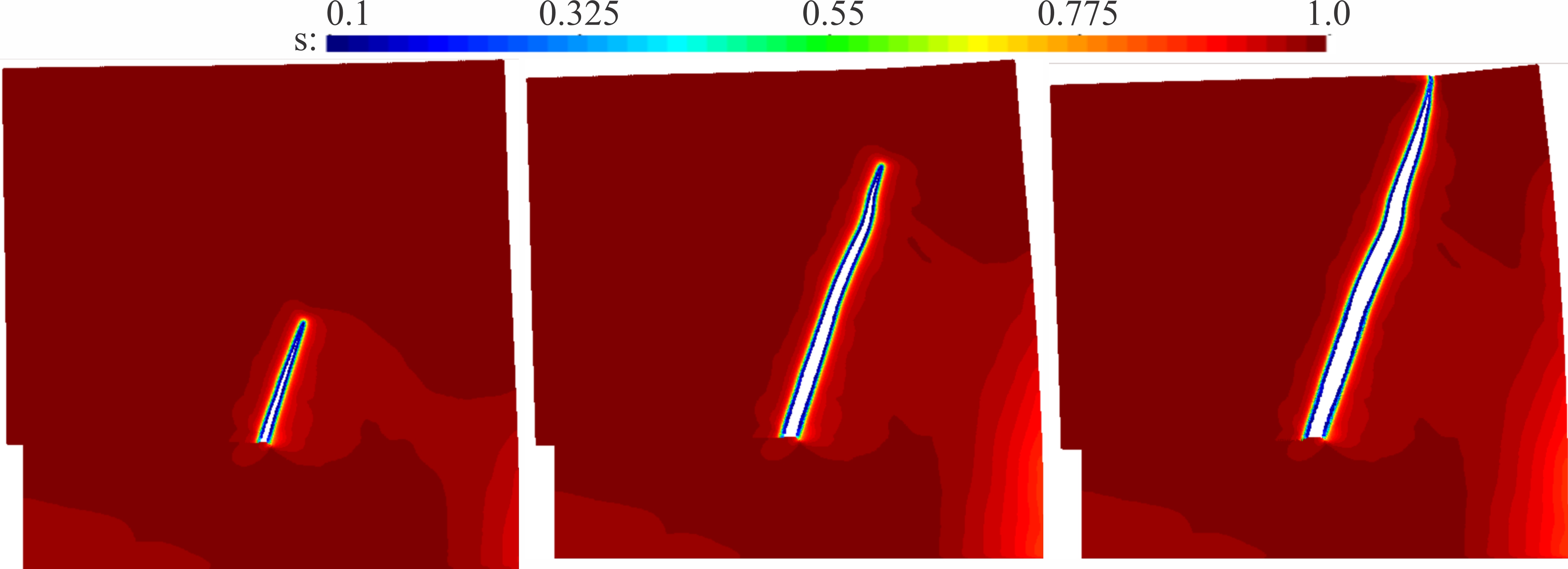}
		\end{center}
		\caption{ }
	\end{subfigure}
	\caption{\centering Kalthoff--Winkler experiment. a) Contours of the phase field parameter and b) opening of the crack at $50$~$\mu$s (left), $80$~$\mu$s (middle), and $101$~$\mu$s (right). In b), particles with $s<s_l=0.1$ are not shown and the deformation is magnified by a factor of $5$.}
	\label{fig:kalthoffwinkPhase}
\end{figure}

\begin{figure}[!htbp]
	\begin{center}
		\includegraphics[width=0.6\linewidth]{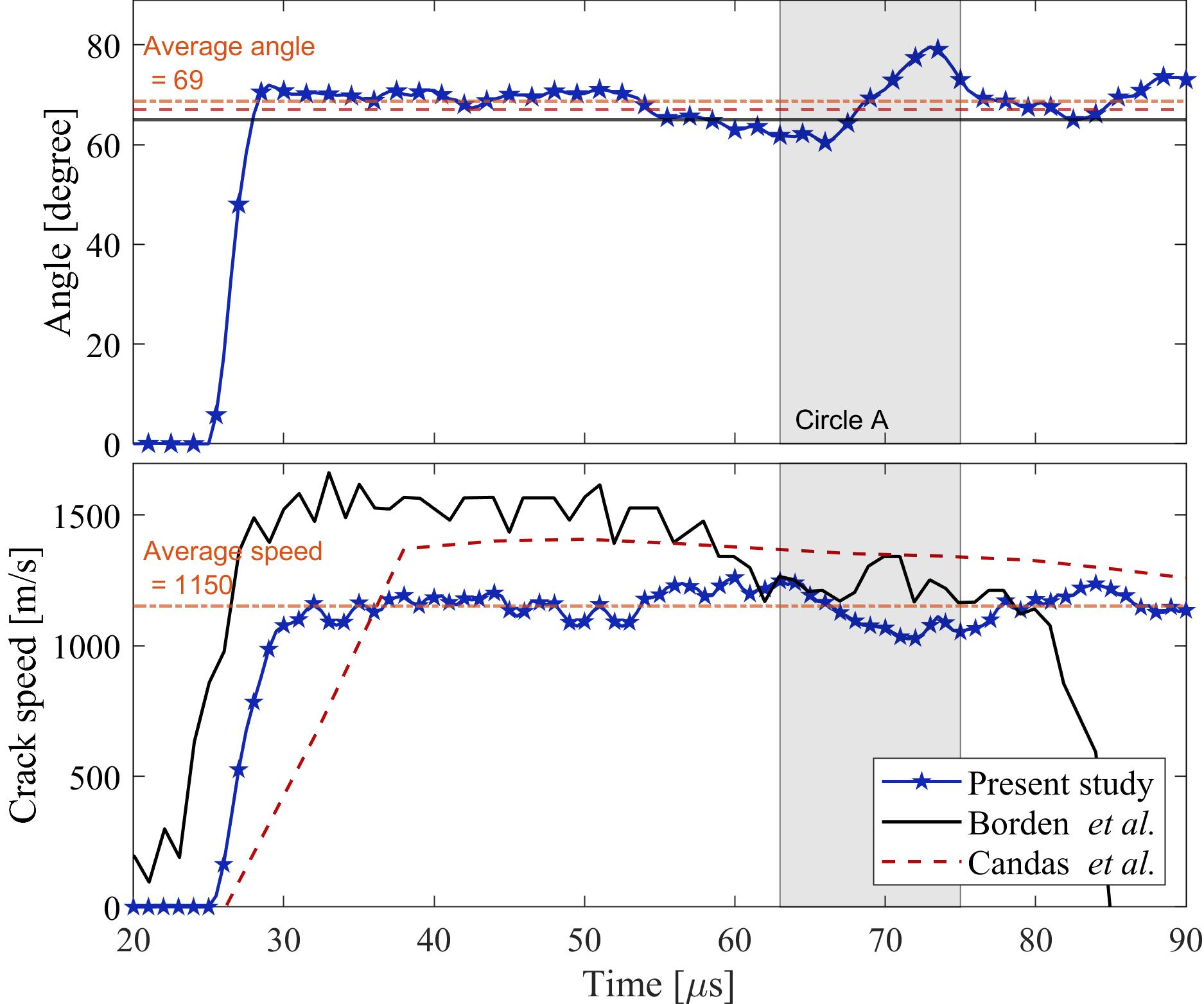}
	\end{center}
	\caption{\centering Kalthoff--Winkler experiment. Propagation speed and orientation of crack compared to phase field results of Borden \etal\cite{Borden2012} and peridynamics results of Candas \etal\cite{Candas2020}.}
	\label{fig:kalthoffwinkGraph}
\end{figure}

Here we revisit the well-known Kalthoff--Winkler experiment where a plate is subjected to an impact shear loading \cite{kalthoff1988failure}. The length and width of the plate are adopted from \cite{Borden2012,Kamensky2018,moutsanidis2018hyperbolic} as $L_x=100$~mm and $L_y=100$~mm, respectively. The material properties are set to $E=190$~GPa, $\rho = 8000$~kg/m$^3$, $\nu=0.3$, $G_c=22.13$~kJ/m$^2$, and the Saint-Venant Kirchhoff constitutive model is employed. The domain is discretized into $638,400$ particles with an initial particle spacing of $\Delta x = 0.125$~mm. The length scale, time step, artificial viscosity parameters, and phase field limit for the soft particles, are chosen as $\epsilon_0=0.125$~mm, $\Delta t=1.0$~ns, $\beta_1=0.04$ and $\beta_2=0$, and $s_l=0.1$, respectively. A geometrical discontinuity with a width of $4\Delta x$ and length of $50$~mm is introduced in the left side of the plate to model the preexisting notch, as shown in \prettyref{fig:kalthoffwinkGeo}.
The plate is subjected to an impact loading in the horizontal direction, applied to the particles located at the lower portion of the left surface. The loading condition is applied as a time-dependent ramp-up velocity of
\begin{equation}
v= \left\{\begin{array}{lr}\frac{t}{t_0} 16.5~m/s & t\leq t_0\\ 16.5~m/s &\text{otherwise}\end{array}\right.\text{ ,}
\end{equation}
in the $X^1$-direction, where $t_0=1$~$\mu$s is the ramp-up time. Upon impact ($t=0$), the displacement field starts propagating towards the right, eventually causing stress and strain energy concentrations at the notch tip, which are then dissipated by the propagating crack (or phase field). \prettyref{fig:kalthoffwinkPhase} shows the snapshots of the results for different stages of the propagation. The crack propagation starts at approximately $25$~$\mu$s and continues diagonally until the full rapture at approximately $100$~$\mu$s. A local disturbance in the diagonal path of the crack is recorded between $t\approx 63$~$\mu$s and $t\approx 75$~$\mu$s (marked in Circle A in \prettyref{fig:kalthoffwinkPhase}) which implies a mixed mode I--II fracture (also reported in \cite{Borden2012,Kamensky2018,moutsanidis2018hyperbolic}). 
\prettyref{fig:kalthoffwinkGraph} compares the propagation speed and crack orientation predicted by the current approach with the ones of the elliptic phase field \cite{Borden2012} and peridynamic \cite{Candas2020} models. Due to the explicit dynamics nature of the present computations that results in a somewhat oscillatory behavior, the graphs were smoothed for visualization purposes. The portion belonging to the disturbance region (Circle A) is marked with gray color.
The predicted results of crack propagation and orientation, as well as the elastodynamic behavior of the plate, are in good agreement with the aforementioned computational efforts and the corresponding experimental work. Overall, the proposed framework is able to capture all the complex qualitative and quantitative characteristics of the crack propagation. 

\subsection{Notched plate with hole}
\begin{figure}[!htbp]
	\begin{center}
		\includegraphics[width=0.35\linewidth]{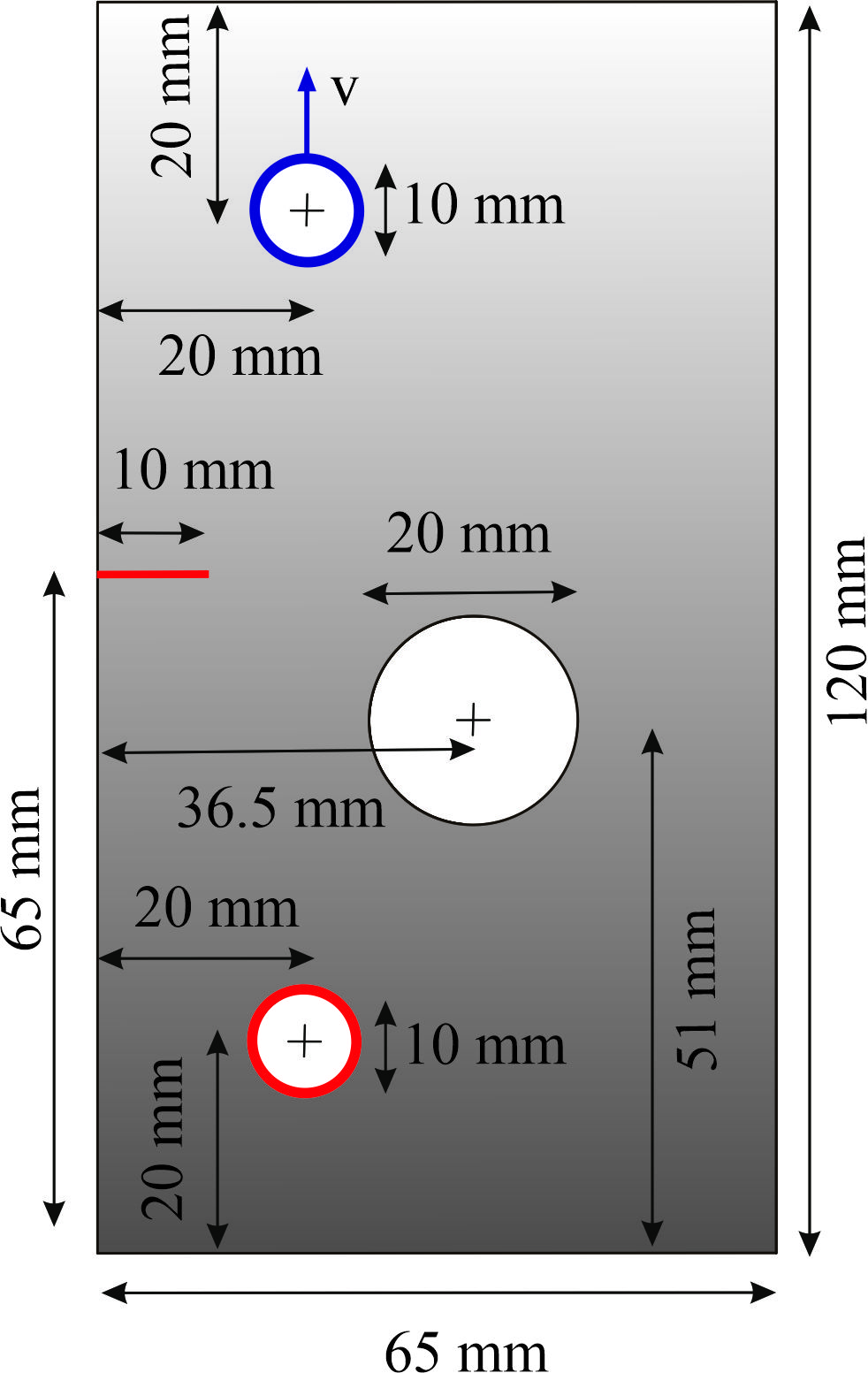}
	\end{center}
	\caption{\centering Notched plate with hole. Geometry and loading conditions. The particles around the lower hole (colored in red) are fixed in the $X^1$- and $X^2$-directions, and the prescribed velocity condition is applied to the particles around the upper hole (colored in blue) in the $X^2-$direction. The preexisting crack with a length of $10$~mm can be seen in red $65$~mm from the bottom edge of the specimen.}
	\label{fig:PlatewithHoleGeo}
\end{figure}

In this example, we simulate the interaction of cracks with other structural defects in a porous plate with an eccentric hole. The geometry of the plate is provided in \prettyref{fig:PlatewithHoleGeo}. The material properties are set to $\lambda=1.94$~GPa, $\mu=2.45$~GPa, $\rho = 2000$~kg/m$^3$, $G_c=2.28$~kJ/m$^2$, and the Saint-Venant Kirchhoff constitutive model is employed. The domain is discretized into $732,852$ particles with an initial particle spacing of $\Delta x = 0.1$~mm. The length scale, time step, artificial viscosity parameters, and phase field limit for the soft particles, are chosen as $\epsilon_0=0.1$~mm, $\Delta t=4.0$~ns, $\beta_1=0.04$ and $\beta_2=0$, and $s_l=0.1$, respectively. 
The original experiment \cite{ambati2015review} was performed in a displacement-control manner, where a prescribed displacement rate of $0.1$~mm/min  was applied. The corresponding simulations \cite{ambati2015review,Navidtehrani2021,Egger2019} were carried out under static displacement loading condition. For both the experiment and the simulations, the prescribed displacement was applied to the upper hole (colored in blue in \prettyref{fig:PlatewithHoleGeo}) of the plate, while the lower hole (colored in red in \prettyref{fig:PlatewithHoleGeo}) was fixed in the horizontal and vertical directions. Here, due to the explicit nature of the SPH method, we conduct dynamic simulations, with three different cases of prescribed velocity in the vertical direction; $v=0.5$~m/s, $0.2$~m/s, $0.1$~m/s, applied on two layers of particles around the upper hole, whereas the two layers of particles around the lower hole are fixed in the $X^1$- and $X^2$-directions. A preexisting crack with a length of $10$~mm is located in the left-central portion of the plate, as shown in \prettyref{fig:PlatewithHoleGeo}, and is modeled by restricting the neighbor search of the particles around the discontinuity region, as suggested in \prettyref{sec:procedure} and \prettyref{fig:cracktypes}d.

\begin{figure}[!htbp]
	\begin{center}
		\includegraphics[width=0.65\linewidth]{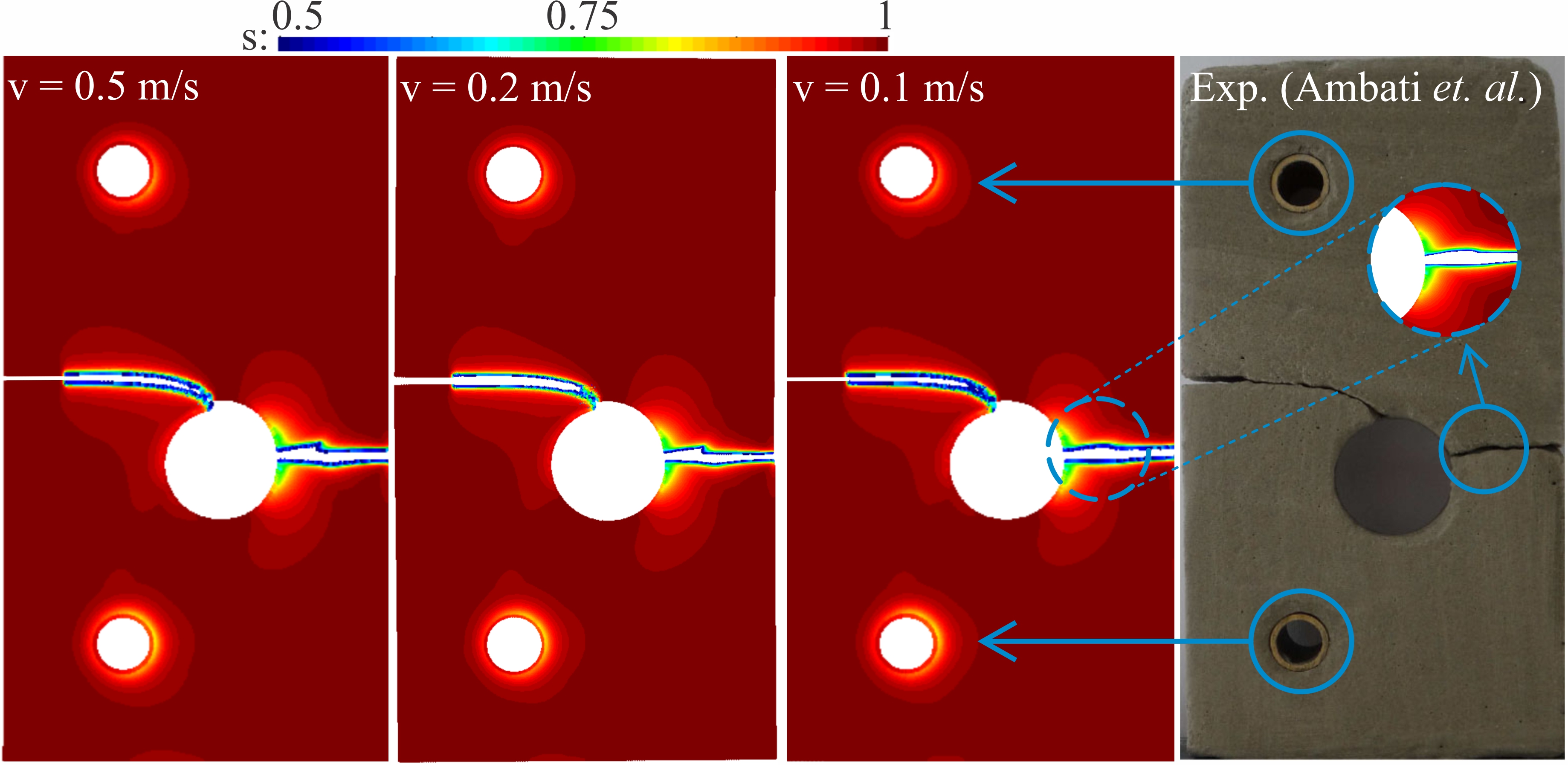}
	\end{center}
	\caption{\centering Notched plate with hole. Comparison of the final crack path with the experimental results of Ambati \etal \cite{ambati2015review}. Particles with $s<0.5$ are not shown.}
	\label{fig:PlatewithHolefinalcrack}
\end{figure}

\begin{figure}[!htbp]
	\begin{center}
		\includegraphics[width=0.65\linewidth]{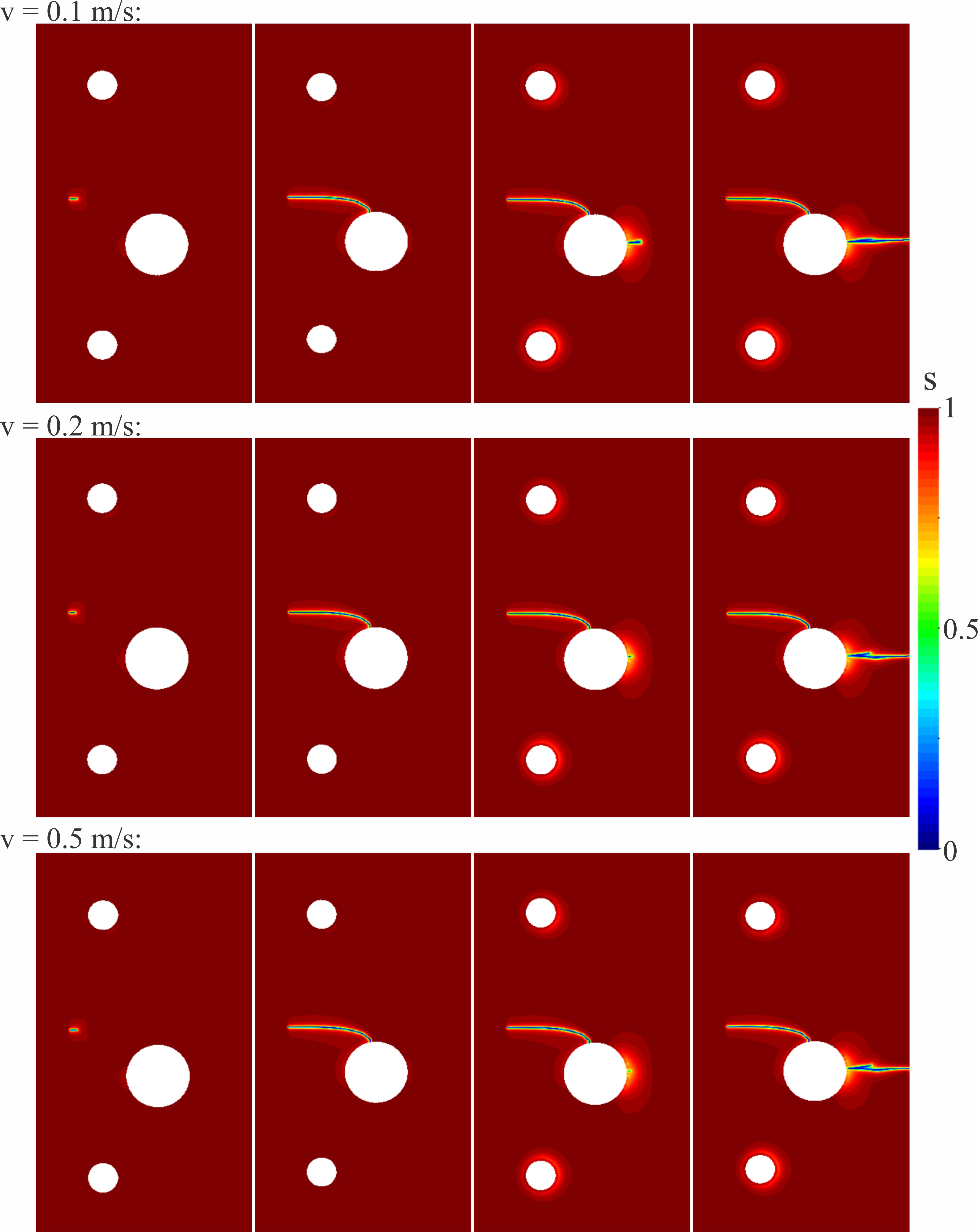}
	\end{center}
	\caption{\centering Notched plate with hole. Contours of the phase field parameter (crack path) under different velocity loading conditions for upper hole displacement values of $0.3$~mm (first column), $0.42$~mm (second column), $1.42$~mm (third column), and $1.48$~mm (fourth column).}
	\label{fig:PlatewithHolevelocitycomp}
\end{figure}

\begin{figure}[!htbp]
	\begin{center}
		\includegraphics[width=0.9\linewidth]{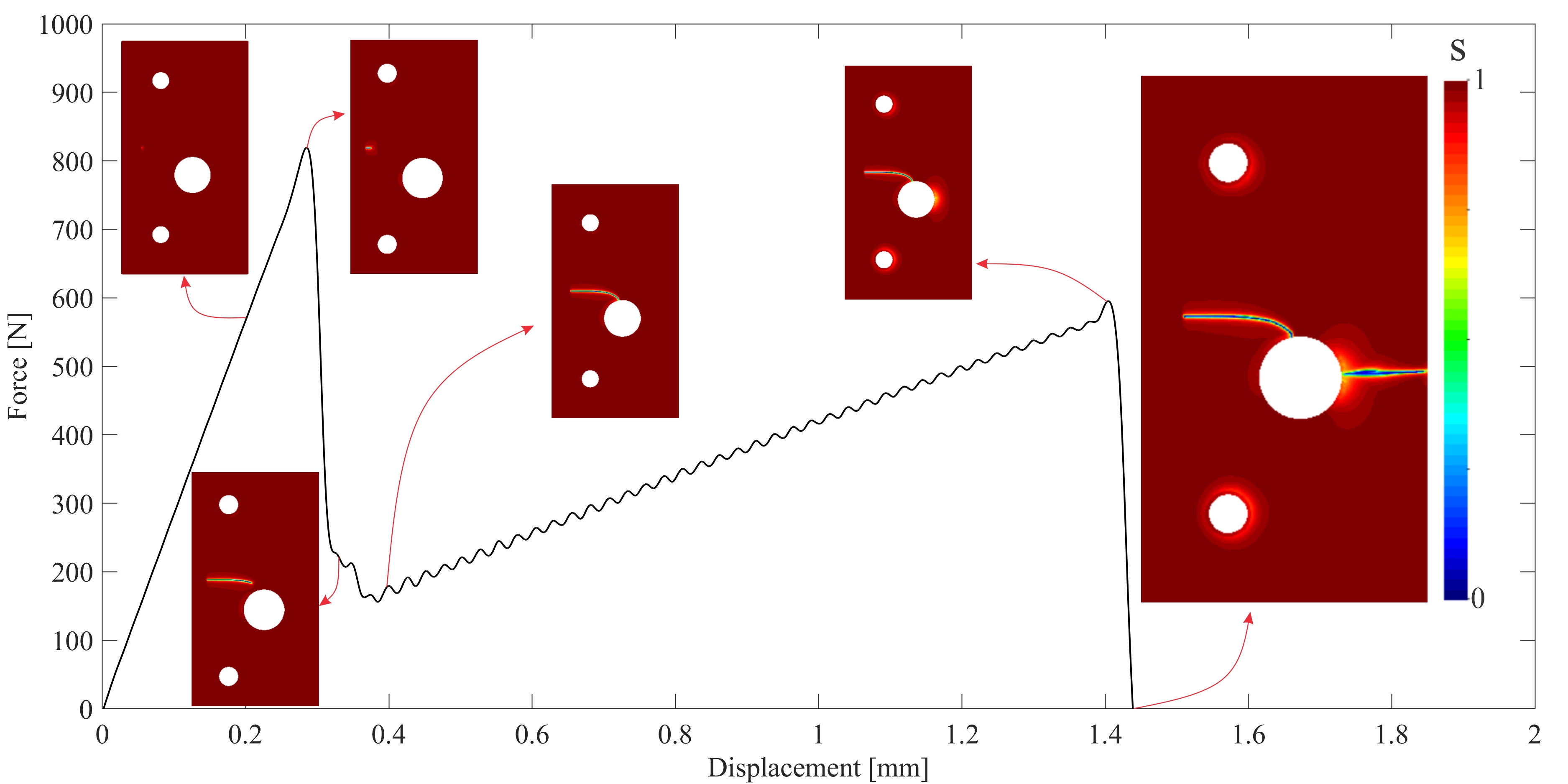}
	\end{center}
	\caption{\centering Notched plate with hole. Reaction force--displacement graph for $v=0.1$~m/s. Snapshots of the crack at various stages are superimposed for a comparative view.}
	\label{fig:PlatewithHolegraph}
\end{figure}

The crack initiation is observed after the upper hole has displaced by approximately $0.3$~mm in the $X^2$-direction for all three cases of the prescribed velocity. The presence of the eccentric hole forms a weak-zone in the vicinity of the crack tip which generates the so-called accelerating effect on the propagating crack \cite{Rahimi2021,Rahimi2020}. As a result, the crack path is directed towards the hole, and eventually joins the hole after a displacement of approximately $0.4$~mm. Upon joining the hole the crack dissipates its energy, leading to a significant delay in propagation, caused by the arresting effect of the hole. At a displacement of approximately $1.4$~mm, a new crack initiates from the right side of the hole and causes a fast and complete rapture. The complex physics in terms of the accelerating and arresting effects of the hole on the crack are well-captured by the proposed approach. Identical accelerating and arresting effects were reported in \cite{Rahimi2021,Rahimi2020} using peridynamics. \prettyref{fig:PlatewithHolefinalcrack} depicts the snapshots of the final crack pattern under different prescribed velocity loading conditions. 
As it is seen, the present approach is able to capture all the major qualitative features of fracture, such as the damaged region around the supports, as well as the tendency of the crack to branch in the vicinity of the eccentric hole (marked with circles in \prettyref{fig:PlatewithHolefinalcrack}), which demonstrates a good agreement with the experimental results of \cite{ambati2015review}. A small deviation of the crack path in terms of small branching is observed for higher prescribed velocity values. This is in agreement with the reported effects of the dynamic loading on crack branching \cite{Bobaru2015} and the experimental observations in \cite{ambati2015review}. \prettyref{fig:PlatewithHolevelocitycomp} compares the phase field values for different loading conditions at various stages of the propagation. As can be seen, for lower values of the prescribed velocity, the branching effect close to the hole vanishes, which agrees well with the findings of \cite{Bobaru2015}.
\prettyref{fig:PlatewithHolegraph} presents the smoothed graph of the reaction force (at the lower hole) versus maximum displacement (displacement of the upper hole) for $v=0.1$~m/s.
We would like to point out that although this problem has been investigated by several researchers in the past \cite{ambati2015review,Navidtehrani2021,Egger2019}, the corresponding reaction force--displacement graphs are not in agreement among each other. Furthermore, all the previous works, in contrast to the present effort, solved the problem in a static manner. Hence, some discrepancies are expected in the reaction force--displacement curve. However, the qualitative features of the graph are captured well. Moreover, the time at which the first crack initiation occurs (sudden drop in the curve at approximately $0.3$~mm) is in excellent agreement with the one reported in \cite{Navidtehrani2021,Egger2019}. 

\subsection{Impact of a spherical projectile on a notched circular plate} \label{sec:impactprobelm}

\begin{figure}[!htbp]
	\begin{center}
		\includegraphics[width=\linewidth]{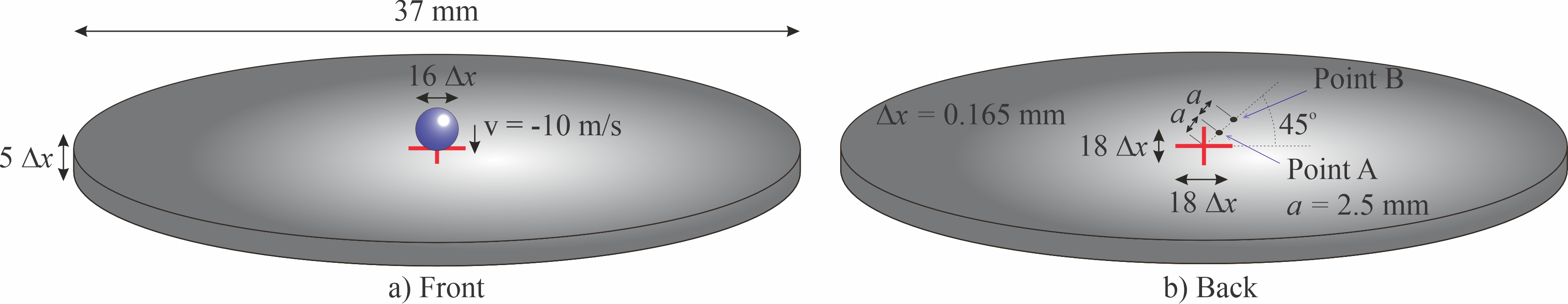}
	\end{center}
	\caption{\centering Impact of a spherical projectile on a notched circular plate. Geometry and loading conditions.}
	\label{fig:diskcontact_geo}
\end{figure}

\begin{figure}[!htbp]
	\begin{center}
		\includegraphics[width=\linewidth]{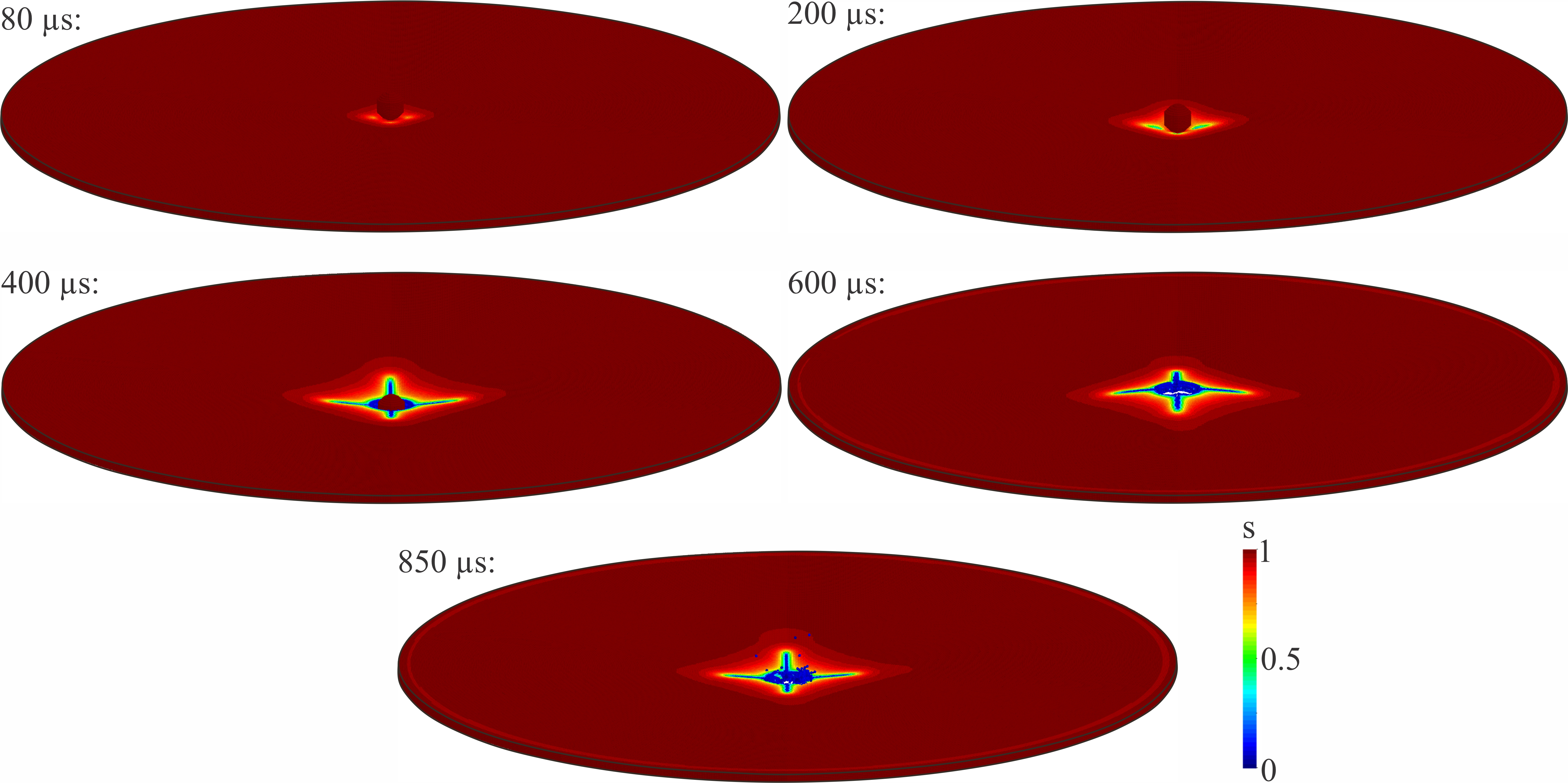}
	\end{center}
	\caption{\centering Impact of a spherical projectile on a notched circular plate. Three dimensional view of the phase field parameter.}
	\label{fig:diskcontact_phasefield_iso}
\end{figure}

\begin{figure}[!htbp]
	\begin{center}
		\includegraphics[width=\linewidth]{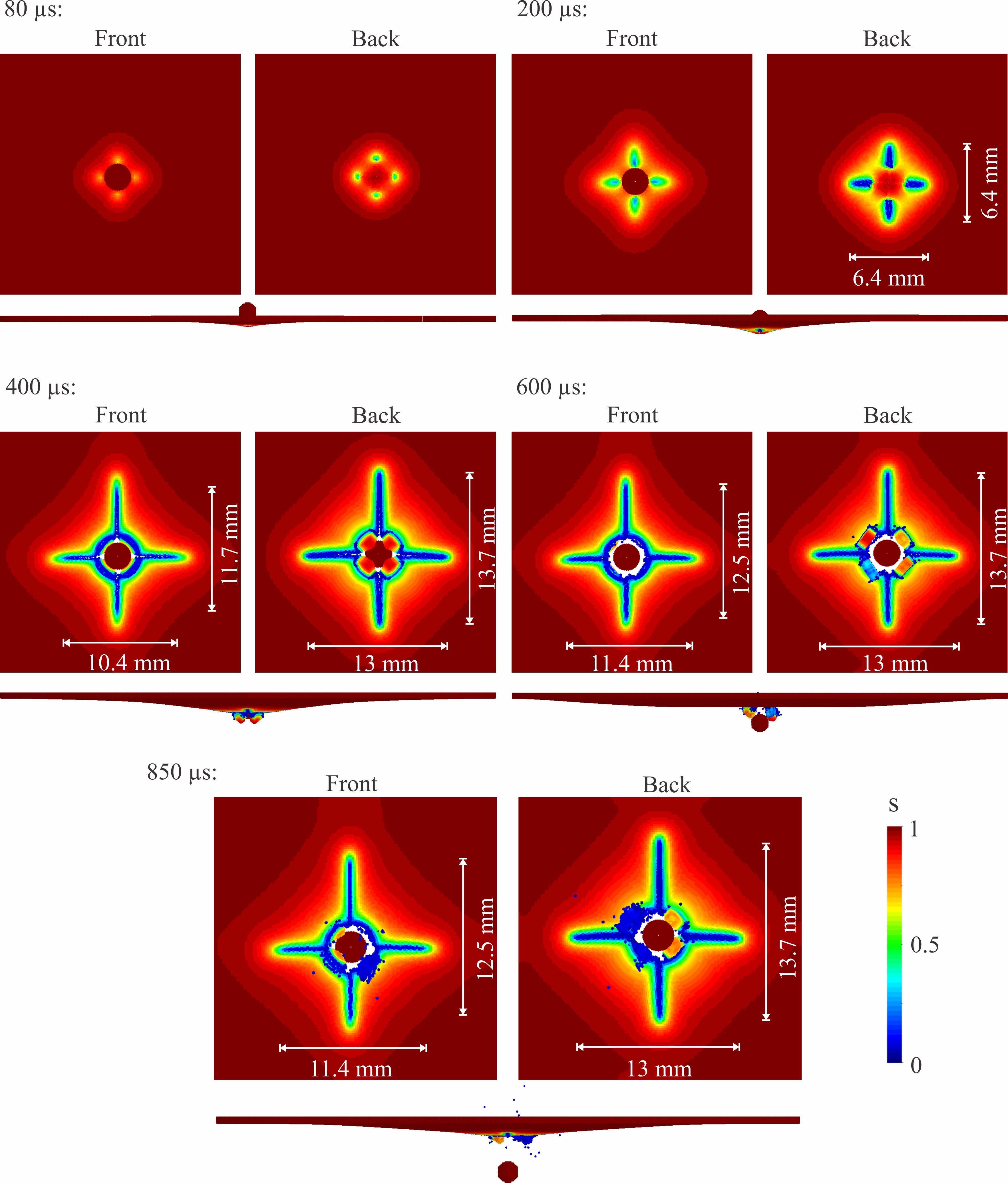}
	\end{center}
	\caption{\centering Impact of a spherical projectile on a notched circular plate. Two dimensional view of the phase field parameter.}
	\label{fig:diskcontact_phasefield_2d}
\end{figure}

\begin{figure}[!htbp]
	\begin{center}
		\includegraphics[width=0.7\linewidth]{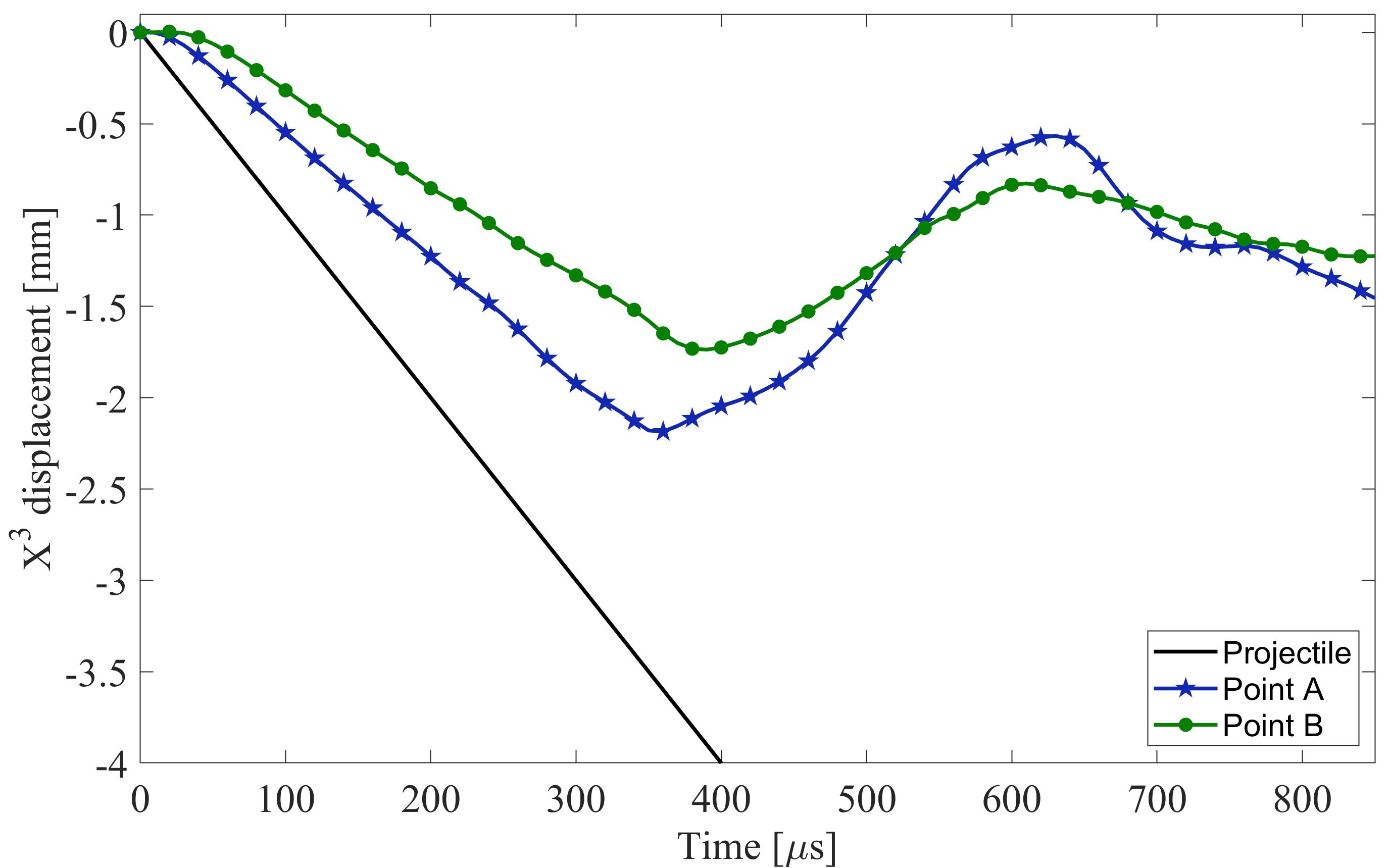}
	\end{center}
	\caption{\centering Impact of a spherical projectile on a notched circular plate. Displacements of Point A and Point B in the $X^3$-direction.}
	\label{fig:diskcontact_graph}
\end{figure}

To further demonstrate the capabilities of our proposed framework we simulate the impact of a rigid spherical projectile on a notched circular plate. \prettyref{fig:diskcontact_geo} illustrates the geometric properties of the domain in which the diameter and thickness of the plate are given as $37$~mm and $5\Delta x$, respectively.
A rigid spherical projectile with a diameter of $16\Delta x$ and a prescribed velocity of $-10$~m/s in the $X^3$-direction impacts the plate at time $t=0$. A through-the-thickness cross-shaped notch with side lengths of $18\Delta x$ is applied at the center of the plate and is modeled by restricting the neighbor search of the particles around the notch. The displacement is constrained at the outer edge of the plate. We set the material and simulation parameters to $\lambda=4.75$~GPa, $\mu=0.9$~GPa, $\rho_0 = 1200$~kg/m$^3$, $G_c=703$~J/m$^2$, $\epsilon_0=\Delta x$, $\Delta t=2.5$~ns, $\beta_1=0.05$, $\beta_2=0$, $s_l=0.1$, $r_0=2\Delta x$, and $K_p=1.2\e{12}$. The domain is discretized into $800,774$ particles and the particle spacing is $\Delta x = 0.165$~mm, whereas the Neo-Hookean constitutive model is employed. \prettyref{fig:diskcontact_phasefield_iso} and \prettyref{fig:diskcontact_phasefield_2d} show the three- and two- dimensional views of the phase field in the deformed plate, respectively.
As can be seen, upon impact the plate starts to deform, and after approximately $80$~$\mu$s the propagation initiates from the back side of the preexisting notch tips. Since the back side of the plate undergoes higher tension than the front, the propagation at the back starts earlier. At later stages of the propagation, the plate is forced to break from the boundaries of the weak regions under high tension which correspond to the regions under the impact of the projectile, thus, leading to an ultimate penetration. After the penetration, the plate bounces back and dissipates its accumulated strain energy marking the end of crack propagation. \prettyref{fig:diskcontact_graph} shows the displacements of Point A and Point B (see \prettyref{fig:diskcontact_geo}) in the $X^3$-direction with respect to time.
As can be seen, the proposed approach is able to efficiently capture the pre- and post-failure behaviors of the complex contact problem with complex fracture patterns. This problem is merely a demonstration of capabilities, and although neither experimental nor numerical results exist for comparison purposes, all the qualitative features that are expected from the underlying physics are produced.

\section{Conclusions}
\label{sec:conclusion}
We presented a robust, accurate, and convergent computational framework for modeling dynamic brittle fracture within the SPH method. The proposed approach employs the well-known SPH approximation technique to discretize a coupled system of a recently developed hyperbolic phase field model of brittle fracture and solid mechanics. The use of SPH allows for the simulation of large deformation problems that potentially involve multi-body interactions and fragmentation, whereas the use of a hyperbolic PDE for the phase field governing equation permits explicit integration in time, without having to solve expensive linear systems associated with elliptic models, or introduce severe time step restrictions associated with parabolic models.

The results of fracture simulations performed with the proposed formulation are in excellent agreement with solutions calculated using standard mesh-based numerical techniques, such as the finite element method and isogeometric analysis. At the same time, they are in very good agreement with results from experimental studies. The last numerical example, although it does not have any corresponding experimental or numerical results, is a demonstration-of-capabilities study, and produces all the qualitative features that are expected from the underlying physics.

This proof-of-concept work presents the core formulation for coupling solid mechanics with phase field for fracture within the SPH framework. Immediate future efforts will focus on large deformation fluid--structure interaction problems, functionally graded materials, as well as parallelization based on graphics processing units (GPU).

\section*{Acknowledgments}
The authors would like to thank Stony Brook Research Computing and Cyberinfrastructure, and the Institute for Advanced Computational Science at Stony Brook University for access to the high-performance SeaWulf computing system, which was made possible by a \$1.4M National Science Foundation grant (\#1531492).

\appendix
\section{Interface force for contact problems in SPH}
\label{app:appendixA}

In a typical contact problem in SPH an additional term arising from the contact of two separate bodies is included in the governing equation. Thus, the resulting momentum equation (\prettyref{eq:MomentumEqTLSPHindex}) becomes
\begin{equation}
	\frac{d v^k_\textbf{i}}{dt} = \sum_{\textbf{j}=1}^{N_{\textbf{i}}} 
	m_{0\textbf{j}}
	\left( 
	\frac{P^{ks}_\textbf{i}}{\rho^2_{0\textbf{i}}} +
	\frac{P^{ks}_\textbf{j}}{\rho^2_{0\textbf{j}}} +
	P^{ks}_{v\textbf{ij}}
	\right)
	~\frac{\partial W_{0\textbf{ij}}}{\partial X_{\textbf{j}}^s} + b_{0\textbf{i}}^k + \sum_{\textbf{a}=1}^{\hat{N}_{\textbf{i}}} f_{\textbf{ai}} r_{\textbf{ai}}^k~.
\end{equation}
In the newly included term, $\hat{N}_{\textbf{i}}$ is the total number of particles located at a separate body within the contact distance of particle $\textbf{i}$, $r_{\textbf{ai}}^k$ is the k component of the relative position vector $\textbf{r}_{\textbf{ai}}=\textbf{x}_{\textbf{a}}-\textbf{x}_{\textbf{i}}$ between particles $\textbf{a}$ and $\textbf{i}$ calculated using the Eulerian coordinates $\textbf{x}$ (see \prettyref{fig:contact}), and $f_{\textbf{ai}}$ is a scalar multiplier. Borrowing mostly from \cite{Yan2022contact} we calculate $f_{\textbf{ai}}$ as
\begin{equation}
    f_{\textbf{ai}} =\frac{m_{0\textbf{a}} K_p}{|\textbf{r}_{\textbf{ai}}|^2} \left\{\begin{array}{lr}B_{\textbf{ai}} R_{\textbf{ai}}& |\textbf{r}_{\textbf{ai}}|\leq r_c\\ \textbf{0} &\text{otherwise}\end{array}\right.\text~{,}
\end{equation}
where
\begin{equation}
    B_{\textbf{ai}} = \frac{1+\frac{|\textbf{r}_{\textbf{ai}}|}{r_c} }{1-e} \left( e^{1+\frac{|\textbf{r}_{\textbf{ai}}|}{r_c}}-1 \right)
    \text~{,}
\end{equation}
\begin{equation}
    R_{\textbf{ai}} =\left\{\begin{array}{lr} \frac{(1+1.5\Tilde{q})(2-\Tilde{q})^3} {8} & \frac{d(|\textbf{r}_{\textbf{ai}}|)}{dt} \leq 0\\ 0 &\text{otherwise}
    \end{array}\right.\text~{,}
\end{equation}
\begin{equation}
    \Tilde{q} = \frac{|\textbf{r}_{\textbf{ai}}|}{1.33\Delta x}\text~{.}
\end{equation}
$m_{0\textbf{a}}$ is the initial mass of particle $\textbf{a}$, $r_c$ is the contact distance, and $K_p$ is a constant determining the potential of the contact. The $K_p$ and $r_c$ constants need to be determined carefully through numerical experiments in order to preserve the physical nature of contact and to avoid any penetration.

\section{Supplementary data}
\label{app:appendixB}
Supplementary material related to this article can be found online at https://doi.org/10.1016/j.cma.2022.115191.

\bibliographystyle{unsrt}
\bibliography{main}

\end{document}